\documentclass[reqno, 10pt]{amsart}

\usepackage{epsf,amssymb,times, graphicx, amscd, verbatim, amsxtra}
\usepackage[all]{xypic}
\usepackage[mathscr]{eucal}

\newtheorem{thm}{Theorem}[section]

\newtheorem{prop}[thm]{Proposition}

\newtheorem{lemma}[thm]{Lemma}

\theoremstyle{definition}
\newtheorem{rmk}[thm]{Remark}

\newcommand{\A}{\mathbb{A}}

\newcommand{\C}{\mathbb{C}}
\newcommand{\R}{\mathbb{R}}

\newcommand{\Q}{\mathbb{Q}}
\renewcommand{\H}{\operatorname{H}}
\newcommand{\CW}{\operatorname{CW}}

\newcommand{\ie}{{\em i.e. }}

\newcommand{\ra}{\rightarrow}

\newcommand{\Ov}{\mathcal{O}_v}

\newcommand{\Sw}{\mathscr{S}}

\newcommand{\Gm}{\mathbb{G}_m}
\newcommand{\D}{D^{\times}}
\newcommand{\JL}{\operatorname{JL}}
\newcommand{\Ind}{\operatorname{Ind}\,}
\newcommand{\nInd}{\operatorname{n-Ind}\,}

\newcommand{\Hom}{\operatorname{Hom}}

\newcommand{\Irr}{\operatorname{Irr}}

\newcommand{\M}{\operatorname{M}_{2\times2}}

\newcommand{\ch}{\operatorname{char}\,}
\newcommand{\Id}{\operatorname{Id}}
\newcommand{\tr}{\operatorname{tr}}
\newcommand{\disc}{\operatorname{disc}}

\newcommand{\GSp}{\operatorname{GSp}}
\newcommand{\GL}{\operatorname{GL}}
\newcommand{\GSO}{\operatorname{GSO}}
\newcommand{\GO}{\operatorname{GO}}
\newcommand{\SO}{\operatorname{SO}}
\newcommand{\OO}{\operatorname{O}}
\newcommand{\Sp}{\operatorname{Sp}}

\newcommand{\Vsigma}{V_{\sigma}}

\newcommand{\ee}{\end{enumerate}}

\newcommand{\bp}{\begin{proof}}
\newcommand{\ep}{\end{proof}}

\begin{document}


\title[local-global non-vanishing of theta lifts]%
{Some local-global non-vanishing results\\ for theta lifts from
orthogonal groups}


\author{Shuichiro Takeda}
\address{Department of Mathematics\\
University of Pennsylvania\\
209 South 33rd St.\\
Philadelphia, PA 19104-6395} \email{stakeda@math.upenn.edu}
\urladdr{http://www.math.upenn.edu/\char126 stakeda}


\date{\today}

\begin{abstract}
We, firstly, improve a theorem of B.\ Roberts which characterizes
non-vanishing of a global theta lift from $\OO(X)$ to $\Sp(n)$ in
terms of non-vanishing of local theta lifts. In particular, we will
remove all the archimedean conditions imposed upon his theorem.
Secondly, following Roberts, we will apply our theorem to theta
lifting of low rank similitude groups. Namely we characterize the
non-vanishing condition of a global theta lift from $\GO(4)$ to
$\GSp(2)$ in our improved setting. Also we consider non-vanishing
conditions of a global theta lift from $\GO(4)$ to $\GSp(1)$ and
explicitly compute the lift when it exists.
\end{abstract}

\maketitle


\section{Introduction}\label{S:intro}


Let $X$ be a symmetric bilinear space over a number field of an even dimension $m$ and $\sigma$ an irreducible cuspidal automorphic representation of $\OO(X,\A_F)$. Assume $\Vsigma$ is a space of automorphic forms giving rise to $\sigma$ . Then we consider the
theta lift $\Theta_n(\Vsigma)$ of $\Vsigma$ to $\Sp(n,\A_F)$ for
$n=\frac{m}{2}$. One of the major questions in the theory of theta
correspondence is to show that the lift $\Theta_n(\Vsigma)$ does not
vanish. B.\ Roberts characterizes the global non-vanishing of the
theta lift in terms of the local counterpart $\theta_n(\sigma_v)$.
To be more precise, in \cite{Rob99-2}, Roberts proves that for
$\frac{m}{2}\leq n$ the theta lift does not vanish if and only if
the local theta lift of each local component $\sigma_v$ does not
vanish under various technical assumptions. In particular he assumes
that the signature of $\OO(X)$ at each real place is of the form
$(2p,2q)$. In this paper, first we will completely remove all the
archimedean assumptions imposed upon his theorem, although we
restrict to the case $\frac{m}{2}= n$. To be precise, we prove

\begin{thm}\label{T:main0}
Let $F$ be any number field, and
$\sigma\cong\otimes\sigma_v$ an irreducible cuspidal automorphic representation
of $\OO(X,\A_F)$ with $\dim X=m$ even. Also let $S_f$ be the finite
set of finite places $v$ such that either $\sigma_v$ is ramified,
$v|2$, or $v$ is ramified in the quadratic extension $F(\sqrt{d})$
of $F$, where $d$ is the discriminant of $X$. Assume:
\begin{enumerate}
\item The (incomplete) standard Langlands $L$-function
$L^S(s, \sigma)$ of $\sigma$ does not vanish at
$s\in\{1,2,\dots,\frac{m}{2}\}$. (A pole is permitted).

\item $\sigma_v$ is tempered for all $v\in S_f$.

\item The local theta lift $\theta_{\frac{m}{2}}(\sigma_v)$ to $\Sp(\frac{m}{2}, F_v)$ exists for all
places $v$.
\end{enumerate}
Then the global theta lift $\Theta_{\frac{m}{2}}(\Vsigma)$ to
$\Sp(\frac{m}{2},\A_F)$ does not vanish.
\end{thm}

Here the temperedness assumption for $v|2$ is due to the lack of
Howe duality principle for even residual characteristic, and thus
quite crucial. The other two are due to the lack of the
corresponding result of \cite{Rob98} for the non-tempered case. In
\cite{Rob99-2}, he assumes that $\pi_v$ is tempered for all
non-archimedean places but here we replace the temperedness
condition by the $L$-function condition. (This is not an improvement
of his theorem, but just another way of stating the theorem,
although this makes a slight difference when we apply it to the
similitude case.) Also in \cite{Rob99-2} he \textit{did not} assume
$\pi_v$ is tempered for archimedean $v$, but in \cite[p.301]{Rob01}
he himself pointed out that this is a mistake and it must be assumed
to be tempered. But in this paper, we will show that, after all,
$\pi_v$ does not have to be tempered for archimedean $v$.

Next we apply this theorem, as Roberts did, to theta lifting for
groups of similitudes. Then we prove the following, which is an
improvement of one of the main theorems of \cite{Rob01}.

\begin{thm}\label{T:main1}
Let $X$ be a symmetric bilinear space of $\dim X=4$ over $F$ and
$\sigma$ an irreducible cuspidal automorphic representation of $\GO(X,\A_F)$.
Assume that $\sigma_v$ is tempered for all $v\in S_f$, where $S_f$
is defined in the same way as in Theorem \ref{T:main0}. Then the
global theta lift $\Theta_2(\Vsigma)$ to $\GSp(2, \A_F)$ does not
vanish if and only if the local theta lift $\theta_2(\sigma_v)$ to
$\GSp(2, F_v)$ does not vanish for all places $v$.
\end{thm}

Notice that the group $\GO(X)$ is disconnected and written as
$\GO(X)\cong\GSO(X)\rtimes \{1,t\}$ for some $t\in\GO(X)$ with
$t^2=1$ which acts on $\GSO(X)$ by conjugation. Each irreducible cuspidal
automorphic representation $\sigma$ of $\GO(X, \A_F)$ is ``extended
from" an irreducible cuspidal automorphic representation $\pi$ of the identity
component $\GSO(X, \A_F)$ in the sense explained in Section
\ref{S:GSO}. Let $d$ be the discriminant of $X$. Roberts in
\cite{Rob01} has shown that for the purpose of similitude theta
lifting, we may assume:

\begin{enumerate}
\item If $d=1$, then an irreducible cuspidal automorphic representation $\pi$ of
$\GSO(X,\A)$ with central character $\chi$ is identified with an irreducible
cuspidal automorphic representation $\tau_1\otimes\tau_2$ of
$\D(\A_F)\times\D(\A_F)$, where $D$ is a quaternion algebra over
$F$, and the central characters of $\tau_1$ and $\tau_2$ are both
$\chi$. In this case, we write $\pi=\pi(\tau_1,\tau_2)$.

\item If $d\neq1$, then an irreducible cuspidal automorphic representation $\tau$ of
$\GSO(X,\A)$ with central character $\chi$ is identified with an irreducible
cuspidal automorphic representation $\tau$ of
${B_E}^{\times}(\A_F)$, where $B$ is a quaternion algebra over the
quadratic extension $E=F(\sqrt{d})$ of $F$, and the central
character of $\tau$ is of the form $\chi\circ N_F^E$. In this case,
we write $\pi=\pi(\tau, \chi)$.
\end{enumerate}

Note that for $\tau_1\otimes\tau_2$ and $\tau$, there are
Jacquet-Langlands lifts $\tau_1^{\JL}\otimes\tau_2^{\JL}$ and
$\tau^{\JL}$ to $\GL(2,\A_F)\times\GL(2,\A_F)$ and $\GL(2,\A_E)$,
respectively. Also for each $\pi$ we can consider the conjugate
$\pi^c$ of $\pi$, whose space $V_{\pi^c}$ of representation is of
the form $\{f\circ c: f\in V_{\pi}\}$ where $f\circ c(g)=f(tgt)$. If
$\pi=\pi(\tau_1,\tau_2)$, then $\pi^c=\pi(\tau_2,\tau_1)$, and if
$\pi=\pi(\tau,\chi)$, then $\pi^c=\pi(\tau^c,\chi)$ where $\tau^c$
is the Galois conjugate of $\tau$. We will prove

\begin{thm}\label{T:main2}
Assume that the global theta lift $\Theta_2(\Vsigma)$ to $\GSp(2,
\A_F)$ does not vanish. Then

\begin{enumerate}
\item If $d=1$ and $\sigma$ is extended from $\pi=\pi(\tau_1,\tau_2)$,
then the theta lift $\Theta_1(\Vsigma)$ to $\GSp(1,\A_F)$ does not
vanish if and only if $\tau_1=\tau_2$. Moreover, if this is the
case, $\Theta_1(\Vsigma)$ is the space of an irreducible cuspidal representation
$\Pi$ such that $\Pi^{\vee}=\tau_1^{\JL}=\tau_2^{\JL}$.

\item If $d\neq1$ and $\sigma$ is extended from $\pi=\pi(\tau,\chi)$,
then the theta lift $\Theta_1(\Vsigma)$ to $\GSp(1,\A_F)$ does not
vanish if and only if $\tau^{\JL}$ is the base change lift of an irreducible
cuspidal automorphic representation $\tau_0$ of $\GL(2,\A_F)$ whose
central character is $\chi_{E/F}\chi$, where $\chi_{E/F}$ is the
quadratic character for $E/F$. Moreover, if this is the case,
$\Theta_1(\Vsigma)$ is the space of an irreducible cuspidal representation $\Pi$
such that $\Pi^{\vee}=\tau_0$.
\end{enumerate}
\end{thm}

In light of those two theorems, one interesting thing to investigate
is, of course, when a given $\pi$ can be extended to $\sigma$ so
that $\Theta_2(\Vsigma)\neq0$. For certain cases, it can be shown
that the answer is ``always". Namely,

\begin{thm}\label{T:main3}\hfill
\begin{enumerate}
\item Assume $\pi$ is generic. Then $\pi$ can be extended to
$\sigma$ so that $\Theta_2(\Vsigma)\neq0$ (without any temperedness
assumption).

\item Assume $\pi$ is such that $\pi^c\ncong\pi$ (but not
necessarily generic). If $\pi$ satisfies the temperedness assumption
as in Theorem \ref{T:main1}, then $\pi$ can be extended to $\sigma$
so that $\Theta_2(\Vsigma)\neq0$.
\end{enumerate}
\end{thm}

We make no claim to the originality for this theorem. The first part
easily follows from a theorem of Howe and Piatetski-Shapiro
\cite{H-PS83}, although they do not explicitly state it in this way.
The second part follows from Theorem \ref{T:main1} together with a
theorem of Roberts in \cite{Rob01}, and he states its tempered
version by using his notion of ``global $L$-packets".\\

Finally, to give arithmetic applications of our lifts, it should be
noted that if $F=\Q$, $d=1$, and $X$ is anisotropic, the lift from
$\GO(X)$ to $\GSp(2)$ is known as the Yoshida lift. (See, say,
\cite{BS91}.) Also if $F=\Q$, $d>0$, and $X$ is anisotropic, then
our lifting from $\GO(X)$ to $\GSp(2)$ together with the
Jacquet-Langlands lift gives a way of constructing a holomorphic
Siegel modular form from a holomorphic
Hilbert modular form.\\

This paper is organized as follows. We first set up our notations in
Section \ref{S:notation}. In Section \ref{S:isometry}, we give a
proof of Theorem \ref{T:main0}. In Section \ref{S:similitude}, we
discuss basics of both global and local theta lifting of groups of
similitudes to the extent necessary for our purpose. In Section
\ref{S:GSO} we review the classification the groups $\GSO(X)$ and
$\GO(X)$ when $\dim X=4$ and their representations. In Section
\ref{S:local}, we will explicitly compute the unramified local theta
lift from $\GO(4)$ to $\GSp(1)$. Then finally, in Section
\ref{S:Main}, we will give the proofs of Theorem \ref{T:main1} and
\ref{T:main2}. In Appendix \ref{S:disconnected}, we give the proof
of Theorem \ref{T:main3}.\\

\begin{center}
Acknowledgements
\end{center}

This work is part of the Ph.D thesis of the author. So the author
would like to thank his advisor Ching-Li Chai for helpful advice.
The author would also like to thank Brooks Roberts for reading the
early draft and making various suggestions, and Herv\'{e} Jacquet
and Stephen Kudla for answering several questions. Thanks are also due to Annegret Paul for explaining her work \cite{Paul}, and Dipendra Prasad for referring the author to his work with R. Schulze-Pillot \cite{PSP}. Indeed, it is his suggestion that Proposition \ref{P:gobal_disconnected} might be proven based on \cite{PSP}.


\section{Notations and Preliminaries}\label{S:notation}


In this paper,  $F$ is a local or global field of $\ch F =0$. If
$E$ is a quadratic extension of $F$, then we denote by $N^E_F$ (or
simply by $N$) the norm map, and by $\chi_{E/F}$ the quadratic
character obtained by local or global class field theory.

We work with smooth representations instead of $K$-finite ones.
Namely if $G$ is a reductive group over a global filed $F$, then by
a (cuspidal) automorphic form we mean a smooth (cuspidal)
automorphic form on $G(\A_F)$ in the sense of \cite[Definition
2.3]{Cogdell04}.

If $\pi$ is a representation of a group, a Harish-Chandra module, or
an automorphic representation, then by $V_{\pi}$ we mean a space which realizes
$\pi$. If $\pi$ is a representation of a real Lie group
then we will denote the space of smooth vectors in $V_{\pi}$ by
$V_{\pi}^{\infty}$ and the space of $K$-finite vectors by
$^{K}V_{\pi}$. If $\pi$ is an admissible representation of a real
Lie group, then we denote the underlying Harish-Chandra module by
$\pi^{\H}$, and thus we have $V_{\pi^{\H}}=\ ^KV_{\pi}$. If $F$ is a
local field, we denote by $\Irr(G(F))$ the collection of equivalence
classes of irreducible smooth admissible representations of $G(F)$.
Also we denote by $\Irr(G(F), \chi)$ the collection of equivalence
classes of irreducible smooth admissible representations of $G(F)$
whose central character is $\chi$. For each $\pi\in\Irr(G(F))$, we
denote the contragredient by $\pi^{\vee}$. If $\pi$ is a
Harish-Chandra module, we denote by $\pi^{\CW}$ the
Casselman-Wallach canonical completion of $\pi$. Thus
$^{K}V_{\pi^{\CW}}=V_{\pi}$. (For the Casselman-Wallach canonical
completion, see \cite{Casselman89}.)

For a finite dimensional vector space $X$ over a global field $F$,
we denote $X\otimes_{F}F_v$ by $X(F_v)$ for each place $v$ and
$X\otimes_{F}\A_F$ by $X(\A_F)$. For a natural number $n$,
$\Sw(X(F_v)^n)$ denotes the space of Schwartz-Bruhat functions. We
set $\Sw(X(\A_F)^n):=\otimes'\Sw(X(F_v)^n)$, where
$\otimes'\Sw(X(F_v)^n)$ is the restricted tensor product over all
places with respect to the characteristic function of $\Ov
x_1+\cdots+\Ov x_m$ for $v$ finite, where $\Ov$ is the ring of
integers of $F_v$ and $x_1,\dots,x_m$ is a fixed basis of $X(F_v)$.

The group $\GSp(n)$ is the algebraic group of symplectic similitudes
of rank $n$ over a field $F$. We realize this group as the group of
$2n\times 2n$ matrices given by $\GSp(n)=\{g\in\GL(2n) :\
^tgJg=\nu(g)J\}$ with
$J=\left(\begin{smallmatrix}0&I_n\\-I_n&0\end{smallmatrix}\right)$,
where $I_n$ is the $n\times n$ identity matrix, and
$\nu:\GSp(n)\ra\Gm$ is the multiplier character. Then the kernel of
$\nu$ is denoted by $\Sp(n)$. When we need to make clear that we are
working with $F$ rational points, we write $\Sp(n,F)$ or $\GSp(n,
F)$, but when there is no danger of confusion, we simply write
$\Sp(n)$ or $\GSp(n)$.

Let $X$ be an even dimensional symmetric bilinear space defined over
a field $F$ of even dimension $m$ equipped with a symmetric bilinear
form. Then we denote by $\GO(X)$ the group of symmetric similitudes
and by $\GSO(X)$ its identity component. If $X$ is defined over a
local or global field $F$ of $\ch F= 0$, then we denote by $\disc
X\in F^{\times}/F^{{\times}^2}$ the discriminant of $X$ when $X$ is
viewed as a quadratic form. We let $\chi_X:F^{\times}\ra\{\pm1\}$ be
the quadratic character of $X$, namely
$\chi_X(a)=(a,(-1)^{\frac{m(m-1)}{2}}\disc X)_F$ for $a\in
F^{\times}$, where $(\ ,\ )_F$ is the Hilbert symbol of $F$.

We denote the (local or global) Weil representation for
$\OO(X)\times \Sp(n)$ by $\omega_{n,X}$ or simply by $\omega$ when
$X$ and $n$ are clear from the context.

If $F$ is an archimedean local field, then the Weil representation
$\omega_{n, X}$ on $\Sw(X(F)^n)$ is a smooth Fr\'{e}chet
representation of the group $\Sp(n)\times\OO(X)$ of moderate growth
in the sense of \cite{Casselman89}. We say that
$\sigma\in\Irr(\OO(X))$ and $\Pi\in\Irr(\Sp(n))$ correspond, or
$\sigma$ corresponds to $\Pi$ if there is a non-zero homomorphism of
Harish-Chandra modules from ${(\omega_{n,X})}^{\H}$ to
$(\Pi\otimes\sigma)^{\H}=\;\Pi^{\H}\otimes\sigma^{\H}$, \ie
$\Hom({(\omega_{n,X})}^{\H}, (\Pi\otimes\sigma)^{\H})\neq 0$, where
$\Hom$ means the set of homomorphisms of Harish-Chandra modules. It
is known that the relation
$\Hom({(\omega_{X,n})}^{\H},(\Pi\otimes\sigma)^{\H})\neq 0$ defines
a graph of bijection between subsets of $\Irr(\Sp(n))$ and
$\Irr(\OO(X))$ up to infinitesimal equivalence (the Howe duality
principle), and in particular if $\sigma$ corresponds to $\Pi$, then
such $\Pi$ is unique up to infinitesimal equivalence, namely
$\Pi^{\H}$ is unique, although $\Pi$ might not be unique. In this
case we write $(\Pi^{\H})^{\CW}=\theta_n(\sigma)$, and we call it
the local theta lift of $\sigma$.

Next assume $F$ is non-archimedean. We say that
$\sigma\in\Irr(\OO(X))$ and $\Pi\in\Irr(\Sp(n))$ correspond, or
$\sigma$ corresponds to $\Pi$ if there is a non-zero
$\Sp(n)\times\OO(X)$ homomorphism from $\omega_{n, X}$ to
$\Pi\otimes\sigma$, \ie $\Hom_{\Sp(n)\times\OO(X)}(\omega_{n,X},
\Pi\otimes\sigma)\neq 0$. If the residue characteristic of $F$ is
odd, it is known that the relation
$\Hom_{\Sp(n)\times\OO(X)}(\omega_{n,X}, \Pi\otimes\sigma)\neq 0$
defines a graph of bijection between subsets of $\Irr(\Sp(n))$ and
$\Irr(\OO(X))$ (the Howe duality principle), and in particular if
$\sigma$ corresponds to $\Pi$, then such $\Pi$ is unique. In this
case we write $\Pi=\theta_n(\sigma)$ and call it the local theta
lift of $\sigma$. If $\sigma$ does not correspond to any
$\Pi\in\Irr(\Sp(n))$, then we say that the theta lift of $\sigma$
vanishes and write $\theta_n(\sigma)=0$. If the residue
characteristic of $F$ is even, then in general it is not known if
the same holds. However, in \cite{Rob98} Roberts has shown, among
other things, that if $\sigma$ is tempered and corresponds to
$\Pi\in\Irr(\Sp(n))$ for $n=\frac{m}{2}$, then $\Pi$ is unique
regardless of the residue characteristic of $F$. So in this case, we
denote $\Pi$ by $\theta_n(\sigma)$ even if the residue
characteristic of $F$ is even.

For an irreducible cuspidal automorphic representation $\sigma$ of $\OO(X,\A_F)$,
we denote by $\Theta_n(V_{\sigma})$ the space of theta lifts of
$\sigma$ to $\Sp(n,\A_F)$, \ie the space generated by the forms of
the form $\theta(f,\phi)$ for $f\in V_{\sigma}$ and
$\phi\in\Sw(X(\A_F)^n)$ which are define by
\[
    \theta(f;\varphi)(g)
    =\int_{\OO(X, F)\backslash \OO(X,\A_F)}\theta(g, h;\varphi)f(h) \, dh
\]
for each $g\in\Sp(n,\A_F)$, where $\theta(g, h;\varphi)$ is the
theta kernel defined by $\theta(g,h;\varphi)=\sum_{x\in
X(F)^n}\omega(g,h)\varphi(x)$.


\section{Theta lifting for isometry groups}\label{S:isometry}


In this section, we will give a proof of Theorem \ref{T:main0},
which is essentially the ingenious argument of \cite{Rob99-2}, which
has its origin in \cite{BS91}. What needs to be done to improve the
theorem of Roberts is to prove the following key technical lemma  in
our improved setting.

\begin{lemma}\label{L:local_zeta}
Let $\sigma\cong\otimes \sigma_v$ satisfy (2) and (3) of Theorem \ref{T:main0}. Also let $S_f$ be as in Theorem \ref{T:main0} and $S_\infty$ the set of infinite places. Set $S=S_f\cup S_\infty$.
For each $k>n=\frac{m}{2}$ and for $v\in S$, the local zeta integral
of the local theta lift $\theta_k(\sigma_v)$ of $\sigma_v$ to
$\Sp(k, F_v)$ can be chosen so that it has a pole at
$s=k-\frac{m}{2}$.  Namely, there exist a matrix coefficient $f_v$
of $\theta_k(\sigma_v)$ and a standard $K$-finite
${\chi_V}_v$-section $\Phi_v$ for $\Sp(k,F_v)$ such that the local
zeta integral $Z(s-\frac{1}{2}, f_v, \Phi_v)$ has a pole at
$s=k-\frac{m}{2}$. (See \cite{Rob99-2} for the notations.)
\end{lemma}

\begin{proof}
The proof is this lemma is given by Roberts in \cite{Rob99-2} under
the assumption that $F$ is totally real and the signature of $X$ at
the real place is of the form $(2p,2q)$. To remove those archimedean
conditions, we need to prove various technical lemmas, which we will
prove in the next subsection. Once those technical lemmas are
proven, this lemma can be shown simply by tracing the proof given by
Roberts.
\end{proof}

Once this lemma is obtained, we are ready to prove Theorem
\ref{T:main0}.

\begin{proof}[Proof of Theorem \ref{T:main0}]
Apply the above lemma to the proof of the main theorem of
\cite[p.146-148]{Rob99-2}. (For the simplicity of a pole of the zeta
integral, use Proposition 1.6 and 1.7 of \cite{Ikeda92}, which does
not assume that $F$ is totally real unlike \cite{Kudla-Rallis94}.)
\end{proof}


\subsection{Some technical lemmas on zeta integrals}


What forced Roberts to impose the conditions on the infinite places
in his theorem in \cite{Rob99-2} is the unavailability of Lemma
\ref{L:local_zeta} for ${\chi_X}_v$ nontrivial and $v$ real, and for
$v$ complex. In this subsection, we give proofs of several technical
lemmas that allow us to prove Lemma \ref{L:local_zeta} in full
generality for the archimedean case. So we assume $F=\R$ or $\C$.
There are basically two technical ingredients we need. The first one
is the theory of the zeta integral for symplectic group at the
archimedean place developed in \cite{Kudla-Rallis90}. There, it is
assumed that $F=\R$ and the character for the zeta integral (which
corresponds to our ${\chi_X}_v$) is trivial. The results in
\cite{Kudla-Rallis90} are used in two places in Roberts' argument in
crucial ways. (See Lemma 7.5 and Theorem 7.8 of \cite{Rob99-2}.) So
we need to extent the results of \cite{Kudla-Rallis90} to full
generality. Although, as mentioned in \cite{Kudla-Rallis90}, there
is no doubt that all the arguments there work for $F=\C$, it seems
that the extension to the case with the non-trivial character (\ie
the sign character) is not completely immediate. Thus first we prove

\begin{prop}\label{P:sign}
All the results in \cite{Kudla-Rallis90} hold even with the presence
of the sign character.
\end{prop}
\begin{proof}
In this proof, all the notations as well as the numberings of
propositions, lemmas, etc are as in \cite{Kudla-Rallis90}. First as
in the trivial character case, the line bundle
$E_{s+\rho}=P\backslash(G\times\C)$ (see (3.3.1) of p.105) where
$p\cdot(g,v)=(pg,\chi(a(g))|a(g)|^{s+\rho}v)$ is trivialized by the
nowhere vanishing section $g\mapsto(g, \chi(a(g))|a(g)|^{s+\rho})$.
So each section $\Phi(g,s)$ is identified with a smooth function on
$\Omega$ via
\[
    \Phi(g,s)\mapsto\chi(a(g))|a(g)|^{-s-\rho}\Phi(g,s).
\]
For each $\Phi$, let us write $\widetilde{\Phi}$ for the
corresponding smooth function in $C^{\infty}(\Omega)$, \ie
\[
    \widetilde{\Phi}(g,s)=\chi(a(g))|a(g)|^{s+\rho}\Phi(g,s).
\]
Then as in p.98, the restriction of $\widetilde{\Phi}|_{K_G}$ to
$K_G$ is left invariant under $P\cap K_G$ and so we may view
$\widetilde{\Phi}|_{K_G}$ as a smooth function on $\Omega\cong(P\cap
K_G)\backslash K_G$.

Now we first have to prove Proposition 3.1.1, which characterizes
convergence of the integral in terms of the order of non-vanishing
on the negligible set. Clearly all the computations until Lemma
3.1.3 (p.101) remain to be true for our case. What needs to be
modified is the equation (3.1.14) in p.101. First notice that
$|\Phi(k,s)|=|\widetilde{\Phi}(k,s)|$ for all $k\in K_G$. Thus the
equation (3.1.14) becomes
\[
    \int_{K}\int_{K}\int_{A^{+}}
    |\widetilde{\Phi}(k(a)i(k_1,k_2),s)|\mu(a)^{-\sigma-\rho_n}\delta(a)da\ dk_1\
    dk_2.
\]
The equation (3.1.16) must be modified as
\[
    f(u)=\widetilde{\Phi}(\phi(u)i(k_1,k_2),s)
\]
as a function on $[0,1]^n$. Then everything else works identically
with the trivial character case. This proves Proposition 3.1.1 and
Corollary 3.1.5.

Now that both Proposition 3.1.1 and Corollary 3.1.5 have been
proven, Theorem 3.2.2. in p.104 can be proven if we prove
Proposition 3.2.1. But the proof  for Proposition 3.2.1. works for
our case because of the trivialization of the line bundle
$E_{s+\rho}$.
\end{proof}

The second ingredient we need in order to remove the archimedean
conditions on the result of \cite{Rob99-2} and thus to obtain Lemma
\ref{L:local_zeta} is the description of a $K$-finite coefficient of
the local lift $\theta_k(\sigma)$ for $\sigma\in\Irr(\OO(X,F))$ with
$\theta_n(\sigma)\neq 0$. Roberts uses the description of a
$K$-finite coefficient of $\theta_k(\sigma)$ in terms of its
Langlands data that are ordered in the standard way. In particular
he uses the fact that, under his assumptions on the signature of
$X$, the first entry of the Langlands data of $\theta_k(\sigma)$ is
the character $\chi_X|\cdot|^{k-n}$. (See the remark preceding
Theorem 7.7 as well as Lemma 2.2 of \cite{Rob99-2}.)

However the real point is not that the first entry of the Langlands
data is $\chi_X|\cdot|^{k-n}$ in the standard order as in
\cite{Rob99-2}, but that the $K$-finite vectors of the
representation $\theta_k(\sigma)$ are realized as the image of the
usual intertwining integral on the representation induced from a
representation of a parabolic subgroup whose Levi factor is
$\GL(1)\times\Sp(k-1)$ in which $\GL(1)$ acts by the character
$\chi_X|\cdot|^{k-n}$. Namely, what we need is

\begin{lemma}\label{L:key_lemma}
Assume that $\sigma$ is an irreducible admissible representation of
$\OO(X,F)$ such that $\theta_n(\sigma)$ exists for
$n=\frac{1}{2}\dim X$. Then for $k>n$, $\theta_k(\sigma)$ is
infinitesimally equivalent to the Langlands quotient of (\ie the
image of the intertwining integral on)
$\Ind_{P}^{\Sp(k-1)}(\chi_X|\cdot|^{k-n}\otimes\sigma_2)$, where $P$
is a (not necessarily the standard choice of) parabolic subgroup
whose Levi factor is $\GL(1)\times\Sp(k-1)$, and $\sigma_2$ is a
(not necessarily tempered) representation of $\Sp(k-1)$. (Of course,
$\chi_X|\cdot|^{k-n}\otimes \sigma_2$ is extended from
$\GL(1)\times\Sp(k-1)$ to $P$ by letting the unipotent radical act
trivially.)
\end{lemma}

This can be shown by combining the following two lemmas, the first
of which is essentially due to Paul \cite{Paul}.

\begin{lemma}\label{L:data}
Assume that $\sigma$ is an irreducible admissible representation of
$\OO(X,F)$ such that $\theta_n(\sigma)$ exists for
$n=\frac{1}{2}\dim X$. Then for $k>n$, $\theta_k(\sigma)$ is
infinitesimally equivalent to the Langlands quotient of (\ie the
image of the intertwining integral on) $\Ind_{P_{k_1\ldots
k_t}}^{\Sp(k)}(\tau_1\otimes\cdots\otimes\tau_t\otimes\tau)$ for
some parabolic $P_{k_1\ldots k_t}$ with $k_1=1$ and
$\tau_1=\chi_X|\cdot|^{k-n}$, where $P_{k_1\ldots k_t}$ is a
parabolic subgroup of $\Sp(k)$ whose Levi factor is isomorphic to
$\GL(k_1)\times\cdots\times\GL(k_t)\times\Sp(k-(k_1+\cdots +k_t))$.
\end{lemma}
\begin{proof}
For $F=\R$, this is just (a part of) Theorem 6.2(1) of \cite{Paul}.
The case for $F=\C$ is identical to the proof of \cite[Theorem
6.2]{Paul} by using Induction Principle and computation of LKT. (See
\cite{Adams95}.)
\end{proof}

\begin{lemma}\label{L:integral}
Assume $\pi$ is an irreducible admissible representation of $\Sp(k)$
which is infinitesimally equivalent to the image of the intertwining
integral on $\Ind_{P_{k_1\ldots
k_t}}^{\Sp(k)}(\tau_1\otimes\cdots\otimes\tau_t\otimes\tau_{t+1})$,
where $P_{k_1\ldots k_t}$ is a parabolic subgroup whose Levi factor
is isomorphic to
$\GL(k_1)\times\cdots\times\GL(k_t)\times\Sp(k-(k_1+\cdots+k_t))$.
Let $\sigma_2$ be an irreducible admissible representation of
$\Sp(k-k_1)$ which is infinitesimally equivalent to the image of the
intertwining integral on $\Ind_{P_{k_2\ldots
k_t}}^{\Sp(k-k_1)}(\tau_2\otimes\cdots\otimes\tau_t\otimes\tau_{t+1})$,
where $P_{k_2\ldots k_t}$ is the obvious subgroup of $P_{k_1\ldots
k_t}$. (Here the convergence of this intertwining integral can be
easily shown.) Then $\pi$ is infinitesimally equivalent to the image
of the absolutely convergent intertwining integral
$J:\Ind_{P_{k_1}}^{\Sp(k)}(\sigma_1\otimes\sigma_2)
\ra\Ind_{\overline{P}_{k_1}}^{\Sp(k)}(\sigma_1\otimes\sigma_2)$
given by
\[
    J(f)(g)=\int_{\overline{U}_{k_1}}f(\bar{u}g)\ d\bar{u},
\]
where $P_{k_1}$ is the parabolic subgroup whose Levi factor is
$\GL(k_1)\times\Sp(k-k_1)$ and whose unipotent radical $U_{k_1}$ is
contained in the unipotent radical of $P_{k_1\ldots k_t}$. (Here we
write $\overline{P}_{k_1}=\ ^tP_{k_1}$ and $\overline{U}_{k_1}=\
^tU_{k_1}$)
\end{lemma}
\begin{proof}
The proof is straight forward and left to the reader.
\end{proof}

Then we can describe a $K$-finite coefficient of $\pi$ as follows.

\begin{lemma}\label{L:coefficient}
By keeping the notations of the above lemma, let $M_{k_1}$ be the
Levi factor of $P_{k_1}$. Also let $H:\Sp(k)\times
\Sp(k)\rightarrow\C$ be a function satisfying the following
properties:
\begin{enumerate}
\item $H(u_1mg_1,\bar{u_2}mg_2)=H(g_1,g_2)$ for $u_1\in U_{k_1},
    \bar{u_2}\in \bar{U}_{k_1}, m\in M_{k_1}$;

\item For any $g_1, g_2\in \Sp(k)$ the function $m\mapsto H(mg_1,g_2)$
    is a coefficient for $\sigma\otimes\delta_{P_{k_1}}^{1/2}$;

\item $H$ is $C^{\infty}$ and $K\times K$-finite on the right,
\end{enumerate}
where $\sigma=\sigma_1\otimes\sigma_2$, $\delta_{P_{k_1}}$ is the
module of $P_{k_1}$, and $K$ is the standard maximal compact
subgroup of $\Sp(k)$. Then the function $f$ defined by
\[
    f(g)=\int_{{M_{k_1}}\backslash G}H(hg,h) dh=\int_{{\overline{U}_{k_1}}\times
    K}H(\bar{u}kg,k)dk d\bar{u}
\]
is absolutely convergent and a $K$-finite coefficient of $\pi$.
(Note that by a coefficient of $\pi$, we mean a finite $\C$ linear
combination of matrix coefficients of $\pi$.)
\end{lemma}

This proposition is essentially due to Jacquet, and the $\GL(k)$
version is stated in (5.5) of \cite{Jacquet79} without a proof.
Since no proof is given there, we will provide a detail proof here
using $\Sp(k)$ as our group. (The proof also works for the
non-archimedean case.) First we need the following lemma. (For this
lemma, the author would like to thank H. Jacquet, who kindly showed
a variant of the proof via a private correspondence
\cite{Jacquet05}.)

\begin{lemma}\label{L:coefficient2}
Let $H$ be as in Lemma \ref{L:coefficient}. Then $H$ is of the form
\[
    H(g_1,g_2)=\sum_i \langle f_i(g_1), f_i'(g_2) \rangle,
\]
for some $f_i\in\Ind_{P_{k_1}}^{G}(\sigma)$ and
$f_i'\in\Ind_{\overline{P}_{k_1}}^{G}(\sigma^{\vee})$, both of which
are $K$-finite.
\end{lemma}
\begin{proof}
Let us simply write $G=\Sp(k)$, $M=M_{k_1}$, $P=P_{k_1}$,
$\overline{P}=\overline{P}_{k_1}$, $U=U_{k_1}$, and
$\overline{U}=\overline{U}_{k_1}$, and also write
$\tau=\sigma\otimes\delta_P^{1/2}$. Let $\mathcal{M}$ be the space
of coefficients of $\tau$. Then the group $M\times M$ acts on
$\mathcal{M}$  by
\[
    (m_1,m_2)\cdot\varphi(m)=\varphi(m_2^{-1}mm_1),\quad \text{for}\ \
    m\in M,\ \ (m_1,m_2)\in M\times M.
\]
Notice that $V_{\tau}\otimes V_{\tau^{\vee}}\cong\mathcal{M}$ as
vector spaces via $(v,w)\mapsto f_{v,w}$. (Here $f_{v,w}$ is defined
by $f_{v,w}(g)=\langle\pi(g)v,w\rangle$, where $\langle,\rangle$ is
the canonical pairing.) So via this isomorphism, $M\times M$ acts on
$V_{\tau}\otimes V_{\tau^{\vee}}$. Moreover it is easy to see that
the representation obtained by this action of $M\times M$ on
$V_{\tau}\otimes V_{\tau^{\vee}}$ is isomorphic to
$\tau|_M\otimes\tau^{\vee}|_M$ as a smooth admissible
representation.

Now define a function $F:K\times K\ra\mathcal{M}$ by
\[
    F(k_1,k_2)(m)=H(mk_1,k_2) \text{ for } (k_1,k_2)\in K\times K.
\]
Then via the isomorphism $V_{\tau}\otimes
V_{\tau^{\vee}}\cong\mathcal{M}$, we have
\[
    F\in\Ind_{(K\cap M)\times(K\cap M)}^{K\times K}(\tau_0\otimes\tau_0^{\vee}),
\]
where $\tau_0=\tau|_{K\cap M}$ and
$\tau_0^{\vee}=\tau^{\vee}|_{K\cap M}$. This can be seen as follows.
First notice that for all $g_1, g_2\in G$ and $m\in M$, we have
$H(mg_1, g_2)=H(m^{-1}mg_1, m^{-1}g_2)=H(g_1, m^{-1}g_2)$. Then if
we write $H(mk_1, k_2)=\sum_{i}\langle\tau(m)v_i, w_i\rangle$, we
have, for $m_1,\ m_2\in M\cap K$,
\begin{align*}
    F(m_1k_1, m_2k_2)(m)&=H(mm_1k_1, m_2k_2)\\
                &=H(m_2^{-1}mm_1k_1, k_2)\\
                &=\sum_{i}\langle\tau(m_2^{-1}mm_1)v_i, w_i\rangle\\
                &=(m_1, m_2)\cdot F(k_1, k_2) (m),
\end{align*}
and so $F(m_1k_1, m_2k_2)=(m_1, m_2)\cdot F(k_1, k_2)$.

Now notice that
\[
    \Ind_{(K\cap M)\times(K\cap M)}^{K\times K}(\tau_0\otimes\tau_0^{\vee})
    \cong\Ind_{K\cap M}^K\tau_0\otimes\Ind_{K\cap M}^K\tau_0^{\vee}.
\]
Thus we have
\[
    F(k_1,k_2)=\sum_{i}f_i(k_1)\otimes\tilde{f}_i(k_2),
\]
for some $f_i\in\Ind_{K\cap M}^K\tau_0$ and $f'_i\in\Ind_{K\cap
M}^K\tau_0^{\vee}$. By viewing $f_i(k_1)\otimes f'_i(k_2)\in
V_{\tau}\otimes V_{\tau^{\vee}}$ as an element in $\mathcal{M}$ via
the isomorphism $V_{\tau}\otimes V_{\tau^{\vee}}\cong\mathcal{M}$,
we have
\[
    F(k_1,k_2)(m)=H(mk_1,k_2)=\sum_{i}\langle\tau(m)f_i(k_1), f'_i(k_2)\rangle.
\]
Since $H$ is $K\times K$-finite on the right, we see that $f_i$'s
and $f'_i$'s are $K$-finite on the right. Now we can extend the
domain of $f_i$'s and $f'_i$'s from $K$ to the all of $G$ by
\begin{align*}
    f_i(u_1m_1k_1)&=\tau(m_1)f_i(k_1)\\
    f'_i(\bar{u}_2m_2k_2)
    &=\tau^{\vee}(m_2)\tilde{f}_i(k_2),
\end{align*}
for $u_1\in U$, $\bar{u}_1\in \overline{U}$ and $m_1, m_2\in M$.
Those are well defined and indeed $f_i\in\Ind_{P}^{G}(\sigma)$ and
$f'_i\in\Ind_{\overline{P}}^{G}(\sigma^{\vee})$, because
$\tau=\sigma\otimes\delta_P^{1/2}$ and
$\tau^{\vee}=\sigma^{\vee}\otimes\delta_{\overline{P}}^{1/2}$, and
also they are $K$-finite.
\end{proof}

Now we are ready to prove Lemma \ref{L:coefficient}.

\begin{proof}[Proof of Lemma \ref{L:coefficient}]
Let us simply write $G=\Sp(k), P_{k_1}=P, U_{k_1}=U,
\overline{U}_{k_1}=\overline{U},$ and $M_{k_1}=M$. Let us also write
$\eta=\Ind_{P}^G(\sigma)$ and
$\eta'=\Ind_{\overline{P}}^G(\sigma^{\vee})$, where $\sigma^{\vee}$
is extended from $M$ to $\overline{P}$ by letting $\overline{U}$ act
trivially. Let $J$ and $J'$ be the intertwining operators for $\eta$
and $\eta'$, respectively, as defined in Lemma \ref{L:integral}.
Then $\pi^{\H}\cong(\eta\slash\ker J)^{\H}$ and
$(\pi^{\vee})^{\H}\cong(\eta'\slash\ker J')^{\H}$, where $\ker J$
and $\ker J'$ are characterized by the property that, for all
$K$-finite $f\in\ker J$ and $f'\in\ker J'$,
\begin{align*}
    J(f)(g) &=\int_{\overline{U}}\langle f(\bar{u}g), \tilde{v} \rangle\ d\bar{u}=0
    \quad\text{for all}\quad\tilde{v}\in\sigma^{\vee},\quad\text{and}\\
    J'(f')(g) &=\int_{U}\langle v, f'(ug) \rangle\ du=0
    \quad\text{for all}\quad v\in\sigma.
\end{align*}
Then if we write $\bar{f}$ and $\bar{f'}$ for the images in
$\eta\slash\ker J$ and $\eta'\slash\ker J'$, respectively, then the
canonical pairing of $\pi$ and $\pi^{\vee}$ is given by
\[
    \langle \bar{f},\bar{f'}\rangle=\int_{M\backslash G} \langle f(h),f'(h)\rangle\ dh,
\]
for $f, f'$ $K$-finite. This can be proven as follows. First of all,
clearly the function $g\mapsto \langle f(g),f'(g)\rangle$ is
$M$-invariant, and so the integral makes sense. Second of all, this
integral absolutely converges, because
\begin{align*}
    \int_{M\backslash G} \langle f(h),f'(h)\rangle\ dh
    &=\int_{\overline{P}\backslash G}
    \int_{\overline{U}}\langle f(\bar{u}k),f'(\bar{u}k)\rangle\ d\bar{u}dk\\
    &=\int_{K} \int_{\overline{U}}\langle f(\bar{u}k),f'(k)\rangle\ d\bar{u}dk,
\end{align*}
where the integral $\int_{\overline{U}}\langle
f(\bar{u}k),f'(k)\rangle\ d\bar{u}$ converges absolutely by Lemma
\ref{L:integral}. And finally, the characterizing property of $\ker
J$ and $\ker J'$ guarantees that the integral is independent of the
choice of the representatives of $\bar{f}$ and $\bar{f'}$. Therefore
a coefficient of $\pi$ is a finite $\C$ linear combination of
functions of the form
\[
    g\mapsto \langle \pi(g)\bar{f},\bar{f'}\rangle
         =\int_{M\backslash G} \langle f(hg),f'(h)\rangle\ dh
         =\int_{K} \int_{\overline{U}}\langle f(\bar{u}kg),f'(k)\rangle\ d\bar{u}dk,
\]
where $f\in\eta$ and $f'\in\eta'$ are $K$-finite.

Now if $H$ is a function satisfying all the three properties, then
by Lemma \ref{L:coefficient2} we have
\[
    H(g_1,g_2)=\sum_i \langle f_i(g_1), f_i'(g_2) \rangle
\]
for some $f_i\in\eta$'s and $f_i'\in\eta'$'s, all of which are
$K$-finite. Thus the lemma follows.
\end{proof}

Once all those lemmas are proven, Lemma \ref{L:local_zeta} for
$F=\R$ can be proven by exactly the same computation as in
\cite{Rob99-2}.

Finally to apply Roberts' argument to the $F=\C$ case, we need the
following.

\begin{lemma}
Assume $F=\C$. The local zeta integral $Z(s-1/2,
f_1\otimes\delta_P^{-1/2}, \Phi_1)$ as in Proposition 8.7 of
\cite{Rob99-2} has a simple pole at $k-n$ for suitably chosen
$\Phi_1$. (See \cite{Rob99-2} for the notations.)
\end{lemma}
\begin{proof}
Let $\phi_1\in C_c^{\infty}(\C^{\times})$ have support in the ball
$B=B(1,\epsilon)\subset\C$ of radius $\epsilon<1$ and center $1$,
and $\phi_1(1)=0$. Let $\Phi_1$ be the section obtained from
$\phi_1$ as in the proof of Proposition 8.7 of \cite{Rob99-2}. Then
the proof is completely analogous to the real case as in
\cite[p.188]{Rob99-2}.
\end{proof}

This lemma is the complex analogue of Proposition 8.7 of
\cite{Rob99-2}. The rest of the proof of Lemma \ref{L:local_zeta} is
identical to the real case.


\section{Theta lifting for similitude groups}\label{S:similitude}


In this section, we will first review the theory of both local and
global theta lifting for groups of similitudes, and then discuss
some relations between the two. Main references for similitude theta
lifting are \cite{Rob01}, and \cite{HST}.

Following \cite{HST}, we extend the Weil representation for
$\Sp(n)\times\OO(X)$ to the group
\[
    R=\{(g,h)\in \GSp(n)\times \GO(X) : \nu(g)\nu(h)=1\}.
\]
(See Section 2 of \cite{HST} for more detail.) We denote this
extended Weil representation, again, by $\omega_{n,X}$ or simply by
$\omega$.

First assume $F$ is non-archimedean. For $\sigma\in\Irr(\GO(X))$ and
$\Pi\in\Irr(\GSp(n))$, we say that $\sigma$ and $\Pi$ correspond, or
$\sigma$ corresponds to $\Pi$ if there is a non-zero $R$
homomorphism from $\omega_{n,X}$ to $\Pi\otimes\sigma$, \ie
$\Hom_{R}(\omega_{n,X}, \Pi\otimes\sigma)\neq 0$. Let
$\GSp(n)^{+}=\{g\in\GSp(n):\nu(g)\in\nu(\GO(X))\}$. If the residue
characteristic of $F$ is odd, it is known that the relation
$\Hom_{R}(\omega_{n,X}, \Pi\otimes\sigma)\neq 0$ defines a graph of
bijection between subsets of $\Irr(\GSp(n)^{+})$ and $\Irr(\GO(X))$
(the Howe duality principle). (This follows from Theorem 4.4 of
\cite{Rob96} together with the multiplicity one theorem of
\cite[Theorem 1.4]{Adler_Prasad}.) Unlike the isometry case, it is
still unknown if the group $\GSp(n)^{+}$ can be replaced by
$\GSp(n)$ even for the odd residual characteristic case, although it
is known to be true for certain cases. (See Theorem 1.8 of
\cite{Rob01}.) For our purpose, however, the following is enough.

\begin{lemma}\label{L:Howe_dual1}
Let $X$ be a four dimensional quadratic form over a non-archimedean
local field $F$ of $\ch F=0$. First assume the residual
characteristic of $F$ is odd. Then if $\sigma\in\Irr(\GO(X))$
corresponds to $\Pi\in\Irr(\GSp(2))$, then $\Pi$ is unique. Next
assume the residual characteristic is even. Then the same holds as
long as $\sigma$ is tempered, and in this case $\Pi$ is also
tempered. Hence if such $\Pi$ exists, we write
$\theta_2(\sigma)=\Pi$, assuming $\sigma$ is tempered if the
residual characteristic is even. If no such $\Pi$ exists, we write
$\theta_2(\sigma)=0$.
\end{lemma}
\begin{proof}
This is just a part of Theorem 1.8 of in \cite{Rob01}.
\end{proof}

Next assume $F$ is archimedean. Then just as in the non-archimedean
case, the extended Weil representation $\omega_{n,X}$ on
$\Sw(X(F)^n)$ is defined, which is a smooth Fr\'{e}chet
representation of the group $R$ of moderate growth in the sense of
\cite{Casselman89}. In particular $(\omega^{\H})^{\CW}=\omega$.
Then for $\sigma\in\Irr(\GO(X))$ and $\Pi\in\Irr(\GSp(n))$, we say
that $\sigma$ and $\Pi$ correspond, or $\sigma$ corresponds to $\Pi$
if there is a non-zero homomorphism of Harish-Chandra modules from
${(\omega_{n,X})}^{\H}$ to $((\Pi\otimes\sigma)|_R)^{\H}$, \ie
$\Hom({(\omega_{n,X})}^{\H}, ((\Pi\otimes\sigma)|_R)^{\H})\neq 0$,
where $\Hom$ means the set of homomorphisms of Harish-Chandra
modules for smooth representations of $R$. Although, just as in the
non-archimedean case, the Howe duality for similitude groups is not
known in full generality, we only need the following for our
purposes.

\begin{lemma}\label{L:Howe_dual2}
Let $X$ be a four dimensional quadratic form over $F=\R$ or $\C$. If
$\sigma\in\Irr(\GO(X))$ corresponds to $\Pi\in\Irr(\GSp(2))$, then
$\Pi$ is unique up to infinitesimal equivalence. In this case, we write
$\theta_2(\sigma)=(\Pi^{\H})^{\CW}$. If no such $\Pi$ exists, we write $\theta_2(\sigma)=0$.
\end{lemma}
\begin{proof}
This is again essentially a part of Theorem 1.8 of in \cite{Rob01}.
In \cite{Rob01}, the signature of $X$ is assumed to be of the form
$(p,q)$ with both $p$ and $q$ even, but this assumption is
unnecessary. Also if $F=\C$, this is obvious because in this case we
have $\GSp(n)^{+}=\GSp(n)$.
\end{proof}

Now let us consider the global case, so assume $F$ is a global
field. Just as the isometry case, define the theta kernel by
$\theta(g,h;\varphi)=\sum_{x\in X(F)^n}\omega(g,h)\varphi(x)$ for
$(g,h)\in R(\A)$ and $\varphi\in\Sw(X(\A_F)^n)$. Then for each
automorphic representation $\sigma$ of $\GO(X,\A_F)$ and for $f\in
V_{\sigma}$, consider the integral
$\theta(f;\varphi)(g)=\int_{\OO(X,F)\backslash
\OO(X,\A_F)}\theta(g,h_1h;\varphi)f(h_1h) \, dh_1$, where
$h\in\GO(X,\A_F)$ is any element such that $\nu(g)\nu(h)=1$. For a
suitable choice of the Haar measure $dh_1$ as in \cite{HK92}, it can
be shown that this integral is absolutely convergent. Also the
invariance property of the measure guarantees that this integral is
independent of the choice of $h$. Now set $\GSp(n, \A_F)^{+}=\{g\in
\GSp(n, \A_F) : \nu(g)\in \nu(\GO(X,\A)\}$. Then $\theta(f;\varphi)$
is a function on $\GSp(n, \A_F)^{+}$ which is left $\GSp(n,F)^{+}$
invariant. We can extend this function to an automorphic form on
$\GSp(n, \A_F)$ by insisting that it is left $\GSp(n, F)$ invariant
and zero outside $\GSp(n, F)\GSp(n, \A_F)^+$. We denote this
automorphic form also by $\theta(f;\varphi)$, whose central
character is $\chi^{-1}\chi_V^n$ where $\chi$ is the central
character of $\sigma$. Then we denote by $\Theta_n(V_{\sigma})$ the
space generated by the automorphic forms $\theta(f;\varphi)$ for all
$f\in V_{\sigma}$ and all $\varphi\in\Sw(V(\A_F)^n)$. If
$\Theta_n(V_{\sigma})$ is in the space of non-zero cusp forms, then
each of the irreducible constituents under right translation provides an irreducible
cuspidal automorphic representation of $\GSp(n,\A_F)$. So let us write $\Theta_n(V_{\sigma})=\bigoplus_i\Pi_i$ where each $\Pi_i$ is an irreducible cuspidal automorphic representation of $\GSp(n,\A_F)$. If we write $\sigma\cong\otimes\sigma_v$
and $\Pi_i\cong\otimes\Pi_{i,v}$, then each $\sigma_v^{\vee}$ corresponds
to $\Pi_{i,v}$.

\begin{rmk}\label{R:twist}
We should mention a certain conventional discrepancy found in the
literature. In \cite{Rob96} and \cite{HK92}, the extended Weil
representation is defined for the group
$R'=\{(g,h)\in\GSp(n)\times\GO(V) : \nu(g)=\nu(h)\}$. On the other
hand in \cite{HST} it is defined for our group $R$. Let us denote
the extended Weil representations of $R'$ by $\omega'$. By direct
computation, it can be shown that $\omega'$ is obtained from
$\omega$ via the isomorphism $R'\rightarrow R$ given by
$(g,h)\mapsto (\nu(g)^{-1}g, h)$. Then for the local case if
$\sigma\in\Irr(\GO(X))$ corresponds to $\Pi\in\Irr(\GSp(n))$ via
$\omega$, then $\pi$ corresponds to $\tilde{\Pi}$ via $\omega'$
where $\tilde{\Pi}$ is defined by
$\tilde{\Pi}(g)=\chi(\nu(g))^{-1}\Pi(g)$ for $\chi$ the central
character of $\Pi$.

The choice of $R$ seems to be completely conventional, but the
reader should be aware that it also affects the global theta lift.
Indeed if we use $R'$, then for the integral defining
$\theta(f;\varphi)(g)$, we have to
choose $h$ to be such that $\nu(g)=\nu(h)$. (Note that the integral
in p.389 of \cite{HST} is not quite correct.) Accordingly, the
central character of $\theta(f;\varphi)$ is $\chi\chi^n_V$, which is
proved in \cite[Lemma 5.1.9]{HK92}.
\end{rmk}

To consider the non-vanishing problem for our similitude case, we
first consider the restriction to the isometry case. If $f$ is an
automorphic form on $\GO(X,\A_F)$, then clearly $f|_{\OO(X,\A_F)}$
is an automorphic form on $\OO(X,\A_F)$. The same thing can be said
to automorphic forms on $\GSp(n, \A)$. If $V$ is a space of
automorphic forms on $\GSp(n,\A_F)$, then we let
$V|_{\Sp(n)}=\{f|_{\Sp(n,\A_F)}:f\in V\}$. Then we have

\begin{lemma}\label{L:comparison}
Let $\sigma$ be an automorphic representation of $\GO(X,\A_F)$. Then
\begin{enumerate}
\item
For $f\in V_{\sigma}$,
$\theta_n(f;\varphi)|_{\Sp(n,\A_k)}=\theta_n(f|_{\OO(X,\A_F)};\varphi)$,
where $\theta_n(f|_{\OO(X,\A_F)};\varphi)$ is the isometry theta
lift of $f|_{\OO(X,\A_F)}$.

\item $\Theta(V_{\sigma})\neq 0$ if and only if $\Theta({V_{\sigma}})|_{\Sp(n)}\neq 0$.

\item $\Theta(V_{\sigma})$ is in the space of cusp forms if and only if $\Theta({V_{\sigma}})|_{\Sp(n)}$ is.
\end{enumerate}
\end{lemma}
\begin{proof}

(1). This is obvious.

(2). Assume $\Theta(V_{\sigma})\neq 0$. Then for some $f\in
V_{\sigma}$, $g\in\GSp(n,\A_F)$, and $\varphi\in\Sw(X(\A_F)^n)$, we
have $\theta_n(f;\varphi)(g)\neq 0$. By definition of
$\theta_n(f;\varphi)$, we may assume $g\in\GSp(n,\A_F)^{+}$. Let
$g_1=g\left(\begin{smallmatrix}I_n&0\\0&\nu(g)^{-1}I_n\end{smallmatrix}\right)\in\Sp(n,\A_F)$.
Then for $h\in\GO(X,\A_F)$ with $\nu(g)\nu(h)=1$, we have
\begin{align*}
    \theta(f;\varphi)(g)
    &=\int_{\OO(X,F)\backslash \OO(X,\A_F)} \left(\sum_{x\in X(F)^n}
    \omega(g,h_1h)\varphi(x)\right)f(h_1h) \, dh_1\\
    &=\int_{\OO(X,F)\backslash \OO(X,\A_F)} \left(\sum_{x\in X(F)^n}
    |\nu(g)|^{\frac{mn}{2}}\omega(g_1,1)\varphi(x\circ h^{-1}{h_1^{-1}})\right)f(h_1h) \,
    dh_1\\
&=|\nu(g)|^{\frac{mn}{2}}\int_{\OO(X,F)\backslash \OO(X,\A_F)}
\left(\sum_{x\in X(F)^n}
    \omega(g_1,h_1)\varphi'(x)\right)f'(h_1) \, dh_1\\
    &=|\nu(g)|^{\frac{mn}{2}}\theta_n(f'|_{\OO(X,\A_F)};\varphi')(g_1),
\end{align*}
where $\varphi'\in\Sw(X(\A_F)^n)$ is given by
$\varphi'(x)=\varphi(x\circ h^{-1})$ and $f'=h\cdot f\in
V_{\sigma}$. (For the action of $\omega(g,h_1h)$, see, for example,
p.261 of \cite{Rob01} together with Remark \ref{R:twist} above.) Now
$\theta_n(f'|_{\OO(X,\A_F)};\varphi')=\theta_n(f';\varphi')|_{\Sp(n,\A_F)}\in
\Theta({V_{\sigma}})|_{\Sp(n)}$. This proves the only if part. The
converse is obvious.

(3). First notice that $\GSp(n)\cong\Sp(n)\rtimes\Gm$, and if
$P\subset\Sp(n)$ is a parabolic  subgroup, then $P\rtimes\Gm$ is a
parabolic subgroup of $\GSp(n)$, and every parabolic subgroup of
$\GSp(n)$ is of this form. Then if $N_P\subset P$ is the unipotent
radical of $P$, then $N_P$ is also the unipotent radical of
$P\rtimes\Gm$. Now assume $\Theta({V_{\sigma}})|_{\Sp(n)}$ is in the
space of cusp forms. Then for each $f\in\Theta({V_{\sigma}})$ and
each $N_P$, we have $\int_{N_P(\A_F)}f(nh)\,dn=0$ for all
$h\in\Sp(n,\A_F)$. We have to show $\int_{N_P(\A_F)}f(ng)\,dn=0$ for
each $g\in\GSp(n,\A_F)$. Let $g_1\in\Sp(n,\A_F)$ be as in (2), and
$f'=(g_1^{-1}g)\cdot f$. Then
$f'|_{\Sp(n,\A_F)}\in\Theta({V_{\sigma}})|_{\Sp(n)}$. So
$\int_{N_P(\A_F)}f'(ng_1)\,dn=0$. But
$\int_{N_P(\A_F)}f'(ng_1)\,dn=\int_{N_P(\A_F)}f(ng_1g_1^{-1}g)\,dn=\int_{N_P(\A_F)}f(ng)\,dn$.
The converse is obvious.

\end{proof}

We should also mention the following, whose proof is elementary and
left to the reader.
\begin{lemma}\label{L:comparison2}
Let $\sigma$ be an automorphic representation of $\GO(X,\A_F)$, and
$\sigma_1$ an irreducible constituent of $\{ f|_{\OO(X,\A_F)} : f\in
V_{\sigma}\}$ as an automorphic representation of $\OO(X,\A_F)$. If
we write $\sigma\cong\otimes\sigma_v$ and
$\sigma_1\cong\otimes{\sigma_1}_v$, then each ${\sigma_1}_v$ is an
irreducible constituent of the restriction $\sigma_v|_{\OO(X, F_v)}$
of $\sigma_v$ to $\OO(X, F_v)$.
\end{lemma}


\section{The groups $\GSO(X)$ and $\GO(X)$ and their representations}\label{S:GSO}


In this section, we briefly review the classification of the groups
$\GSO(X)$ and $\GO(X)$ and their representations when $\dim X=4$.
All the proofs are found in \cite{Rob01}, \cite{HST}, or citations
therein.


\subsection{The local case}


First consider the case where $F$ is a local field. Then we have
\begin{lemma}
For the sake of similitude theta correspondence, we may assume
\begin{enumerate}
\item If $d=1$, then $\GSO(X)\cong\GSO(M_2)$ or $\GSO(D)$, where $M_2$ is the space of
$2\times 2$ matrices over $F$ with the quadratic form given by
$-\det$, and $D$ is the unique division quaternion algebra over $F$
made into a quadratic form in the usual way. Then there is a natural
bijection between $\Irr(\GSO(X),\chi)$ and the set of irreducible
admissible representations $\tau_1\otimes\tau_2$ of
$\GL(2)\times\GL(2)$ if $\GSO(X)\cong\GSO(M_2)$, or of
$D^{\times}\times D^{\times}$ if $\GSO(X)\cong\GSO(D)$, such that
both $\tau_1$ and $\tau_2$ have the same central character $\chi$.
In this case, we write $\pi=\pi(\tau_1,\tau_2)$.

\item If $d\neq 1$, then $\GSO(X)\cong\GSO(X_E)$, where $X_E=\{x\in\M(E)| \ ^cx^t=x\}$
is the space of Hermitian matrices over $E=F(\sqrt{d})$ with the
quadratic form given by $-\det$. Then there is a natural bijection
between $\Irr(\GSO(X),\chi)$ and $\Irr(\GL(2,E),\chi\circ N^E_F)$.
In this case, we write $\pi=\pi(\tau,\chi)$.
\end{enumerate}
\end{lemma}

Let $\pi\in\Irr(\GSO(X))$. We define $\pi^c$ by taking
$V_{\pi^c}=V_{\pi}$ and by letting $\pi^c(g)f=\pi(tgt)f$ for all
$g\in \GSO(X)$ and $f\in V_{\pi}$. Also notice that
$\GO(X)\cong\GSO(X)\rtimes\{1,t\}$, where we choose $t$ to act on
$X$ as the matrix transpose if $X=M_2$ or $X=X_E$ and the quaternion
conjugation if $X=D$. Then we have
\begin{itemize}
\item If $\pi\ncong\pi^c$, then $\Ind_{\GSO(X)}^{\GO(X)}\pi$ is irreducible,
and we denote it by $\pi^{+}$.

\item If $\pi\cong\pi^c$, then $\Ind_{\GSO(X)}^{\GO(X)}\pi$ is reducible. Indeed,
it is the sum of two irreducible representations, and we write
$\Ind_{\GSO(X)}^{\GO(X)}\pi\cong\pi^{+}\oplus\pi^{-}$. Here, $t$
acts on $\pi^{\pm}$ via a linear operator $\theta^{\pm}$ with the
property that $(\theta^{\pm})^2=\Id$ and $\theta^{\pm}\circ
g=tgt\circ\theta^{\pm}$ for all $g\in \GSO(X)$.
\end{itemize}

We can be more explicit about the irreducible components $\pi^{+}$
and $\pi^{-}$. First assume $d=1$. In this case, it is easy to see
that, via $\rho$, $t$ acts on $\GL(2)\times\GL(2)$ or
$D^{\times}\times D^{\times}$ by $t\cdot(g_1,g_2)=(g_2,g_1)$, and if
$\pi=\pi(\tau_1,\tau_2)$ is such that $\pi\cong\pi^c$, then
$\tau_1\cong\tau_2$. Then we can choose $\theta^{\pm}$ to be such
that $\theta^{+}(x_1\otimes x_2)=x_2\otimes x_1$ and
$\theta^{-}(x_1\otimes x_2)=-x_2\otimes x_1$ for $x_1\otimes x_2\in
V_{\tau_1\otimes\tau_2}$. We choose $\tau^{+}$ and $\tau^{-}$
accordingly. Note that our choice of $\pi^+$ is the canonical extension defined in
\cite[p.10]{PSP}.

Next assume $d\neq1$. In this case $t$ acts, via
$\rho$, on $\GL(2,E)$ in such a way that $t\cdot g=\ ^cg$. If
$\pi=\pi(\tau,\chi)$ is such that $\pi\cong\pi^c$, then
$\tau\cong\tau^c$. Note that $\tau$ has a unique Whittaker model,
namely it is realized as a space of functions $f:GL(2,E)\ra\C$ such
that
$f\left(\left(\begin{smallmatrix}1&a\\0&1\end{smallmatrix}\right)\right)=\psi_v(\tr\,
a)f(g)$ for all $a\in E$ and $g\in\GL(2,E)$, where $\psi_v$ is a
fixed additive character of $F$. Then we define $\theta^{\pm}$ to be
the linear operator that acts on this space of Whittaker functions
by $f\mapsto \pm f\circ c$, and $\theta^{+}$ is chosen to be the one
that acts as $f\mapsto f\circ c$. We choose $\pi^{+}$ and $\pi^{-}$
accordingly.

\begin{rmk}
We should note that our choice of $\pi^{+}$ and $\pi^{-}$ is
different from that of Roberts in \cite{Rob01}, but rather we follow
\cite{HST}, although, if $\pi$ is spherical, it turns out that our
$\pi^{+}$ is spherical and coincides with the notation of
\cite{Rob01}. Also the reader should notice that in the above
discussion the fields $F$ and $E$ do not have to be non-archimedean.
\end{rmk}


\subsection{The global case}


Now we consider the global case, and hence in this subsection $F$ is
a global field of $\ch F =0$ and the groups $\GSO(X)$, $\GO(X)$,
etc denote algebraic groups over $F$. If $d\neq 1$ and
$E=F(\sqrt{d})$, then let $c$ be the non-trivial element in
$\text{Gal}(E/F)$. For each quaternion algebra $D$ over $F$, let
$B_{D,E}=D\otimes E$. Then for each $g\in B_{D,E}$, we define
$^cg^{\ast}$ by linearly extending the operation $^c(x\otimes
a)^{\ast}= x^{\ast}\otimes\, ^ca$ where $^c$ is the Galois
conjugation and $^{\ast}$ is the quaternion conjugation. Then the
space $X_{D,E}=\{g\in B_{D,E} :\, ^cg^{\ast}=g\}$ can be made into a
four dimensional quadratic space over $F$ via the reduced norm of
the quaternion $B_{D,E}$. Similarly to the local case, we have

\begin{lemma}
For the sake of similitude theta correspondence, we may assume
\begin{enumerate}
\item If $d=1$, then $\GO(X)$ is isomorphic to $\GO(D)$ for some (possibly split)
quaternion algebra over $F$. Then there is a natural bijective
correspondence between an irreducible cuspidal automorphic representation $\pi$
of $\GSO(X,\A_F)\cong \GSO(D, \A_F)$ whose central character is
$\chi$ and an irreducible cuspidal automorphic representation
$\tau_1\otimes\tau_2$ of $\D (\A_F)\times \D(\A_F)$ such that both
$\tau_1$ and $\tau_2$ have the central character $\chi$. In this
case, we write $\pi=\pi(\tau_1,\tau_2)$.

\item If $d\neq 1$, then there exists a quaternion algebra $D$ over $F$ such
that $\GO(X)\cong \GO(X_{D,E})$. Then there is a natural bijective
correspondence between an irreducible cuspidal automorphic representation $\pi$
of $\GSO(X, \A_F) \cong\GSO(X_{D,E}, \A_F)$ whose central character
is $\chi$ and an irreducible cuspidal automorphic representation $\tau$ of
$B_{D,E}^{\times}(\A_E)$ whose central character is of the form
$\chi\circ N^E_F$. In this case, we write $\pi=\pi(\tau,\chi)$.
\end{enumerate}
\end{lemma}

Finally, we consider the relation between irreducible cuspidal automorphic
representations of the two groups $\GSO(X,\A_F)$ and $\GO(X,\A_F)$.
First define $\pi^c$ by taking $V_{\pi^c} = \{f\circ c : f\in
V_{\pi}\}$, where $c:\GSO(X, \A_F)\ra \GSO(X, \A_F)$ is the
isomorphism given by conjugation $g\mapsto tgt$. Then clearly
$\pi^c$ is an irreducible cuspidal automorphic representation of $\GSO(X, \A_F)$.
(Note that as an admissible representation, $\pi^c$ is isomorphic to
the representation $\pi'$ with $V_{\pi'}=V_{\pi^c}$ and the action
defined by $\pi'(g)f=\pi(tgt)f$, and so if we write
$\pi\cong\otimes\pi_v$, then $\pi^c\cong\otimes\pi_v^c$.) By
multiplicity one theorem, $\pi\cong\pi^c$ implies
$V_{\pi}=V_{\pi^c}$ and in this case $f\circ c\in V_{\pi}$. Also let
$\sigma$ be an irreducible cuspidal automorphic representation of $\GO(X, \A_F)$.
Define $V_{\sigma}^{\circ}=\{f|_{\GSO(X, \A_F)} : f\in
V_{\sigma}\}$. Then either $V_{\sigma}^{\circ}=V_{\pi}$ for some irreducible
cuspidal automorphic representation $\pi$ of $\GSO(X, \A_F)$ such
that $\pi=\pi^c$, or $V_{\sigma}^{\circ}=V_{\pi}\oplus V_{\pi^c}$
for some irreducible cuspidal automorphic representation $\pi$ of $\GSO(X,
\A_F)$ such that $\pi\neq\pi^c$. (See \cite[p.381--382]{HST}.) Then we have

\begin{prop}\label{P:gobal_disconnected}
Define $\widehat{\pi}$ to be the sum of all the irreducible cuspidal automorphic
representations of $\GO(X, \A_F)$ lying above $\pi$, \ie
$\widehat{\pi}=\oplus_i\sigma_i$ where $\sigma_i$ runs over all the irreducible
cuspidal automorphic representations of $\GO(X, \A_F)$ such that
$V_{\sigma_i}^{\circ}=V_{\pi}$ if $\pi=\pi^c$, or $V_{\sigma_i}^{\circ}=V_{\pi}\oplus V_{\pi^c}$ otherwise.

Then
\[
    \widehat{\pi}\cong\bigoplus_{\delta}\underset{v}{\otimes}\pi_v^{\delta(v)},
\]
where $\delta$ runs over all the maps from the set of all places of
$F$ to $\{\pm\}$ with the property that $\delta(v)=+$ for almost all
places of $F$, and $\delta(v)=+$ if $\pi_v\ncong\pi_v^c$, and
further if $\pi\cong\pi^c$, then $\underset{v}{\prod}\delta(v)=+$.
Moreover each $\otimes\pi_v^{\delta(v)}$ is (isomorphic to) an irreducible
cuspidal automorphic representation of $\GO(X,\A_F)$.
\end{prop}
\begin{proof}
When $\pi$ is generic, the proof is given in \cite[p.382-383]{HST}. Also one can see that the proof there works even if $\pi$ is not generic as long as $\pi$ satisfies $\pi\ncong\pi^c$. So we assume that $\pi$ is not generic and $\pi\cong\pi^c$. The basic idea of the proof is in \cite[p.10-11]{PSP}. Note that our
choice of $\pi^+$ is the canonical extension defined in
\cite[p.10]{PSP}. Then following \cite{PSP}, define a linear
functional $e_\pi:V_\pi\rightarrow\C$ given by evaluation at $1$,
\ie $e_\pi(f)=f(1)$. Now notice that we have an isomorphism
$V_\pi\cong V_{\pi^\delta}$ of vector spaces, because
$\pi\cong\otimes\pi_v$ and for each place $v$, $V_{\pi_v}\cong
V_{\pi_v^{\delta(v)}}$. Then via this isomorphism we define a linear
functional $e_{\pi^\delta}:\pi^\delta\rightarrow\C$. Note that there
is a natural action of $\GO(X, \A_F)$ on $e_{\pi^\delta}$ via
$g\cdot e_{\pi^\delta}(f)=e_{\pi^\delta}(g^{-1}\cdot f)$. Then one
can see that $\pi^\delta$ is automorphic if and only if $g\cdot
e_{\pi^\delta}=e_{\pi^\delta}$ for all $g\in\GO(X,F)$. But clearly
for all $g\in\GSO(X,F)$, $g\cdot e_{\pi^\delta}=e_{\pi^\delta}$.
Hence it suffices to show that $t\cdot
e_{\pi^\delta}=e_{\pi^\delta}$. (Note that that
$t=\otimes_v t_v$ is coherent in the sense of \cite[Middle of
p.11]{PSP}. Also note that in \cite{PSP} our $t_v$ is denoted by $\sigma_v$.) Then it can be easily seen that $t_v\cdot e_{\pi^\delta}=\delta(v)e_{\pi^\delta}$. Then $t=\otimes_v t_v$ acts
as multiplication by $\prod_v\delta(v)$. Thus $e_{\pi^\delta}$ is
invariant under the action of $t$ if and only if
$\prod_v\delta(v)=1$. This completes the proof.
\end{proof}

This proposition tells us that if $\pi$ is an irreducible cuspidal automorphic
representation of $\GSO(X, \A_F)$ and  $\delta$ is a map from the set of all places of
$F$ to $\{\pm\}$ having the property described in the above
proposition, then there is an irreducible cuspidal automorphic representation
$\sigma=(\pi, \delta)$ of $\GO(X, \A_F)$ lying above $\pi$ such that
$\sigma\cong\otimes\pi_v^{\delta(v)}$. We call such a map $\delta$ an
``extension index'' of $\pi$, and $(\pi, \delta)$ the extension of
$\pi$ with an extension index $\delta$.


\section{Local parameters of unramified theta lifts}\label{S:local}


After some preliminaries, we will explicitly compute the local
parameters of the unramified theta lifts from $\GO(4)$ to
$\GSp(1)(=\GL(2))$.

In this section, the groups $\GO(X,F_v), \GSp(n, F_v)$, etc are all
denoted simply by $\GO(X)$, $\GSp(n)$, etc, and we simply write $F$
for $F_v$. Moreover we assume that $v$ is finite. Also ``$\Ind$''
always means unnormalized induction, and whenever we use normalized
induction, we use the notation ``$\nInd$".


\subsection{Preliminaries}


For our computation of the local parameters, we need the Jacquet
module of the Weil representation, which is done in \cite{HST},
which, in turn, comes from \cite{Kudla86}. We will repeat the
essential point. For this, let us decompose $X$ as $X=Y_r\oplus W
\oplus Y_r^{\ast}$, where $Y_r$ is a totally isotropic space and
$Y_r^{\ast}$ is its complement so that $Y_r\oplus Y_r^{\ast}$ is $r$
copies of the hyperbolic plane. We denote the standard basis of
$Y_r$ by $\{f_1,f_2,\dots, f_r\}$, and write $l=\dim W$ so that
$m=2r+l$. Now let $Q_r$ be the parabolic subgroup of $\GO(X)$
preserving the flag $\langle f_1 \rangle\subset \langle f_1,f_2
\rangle\subset\dots\subset \langle f_1,f_2,\dots,f_r \rangle$ so
that its Levi factor is isomorphic to $\GO(W)\times\Gm^r$. Let $Q$
be the parabolic subgroup preserving the flag $\langle
f_1,f_2,\dots,f_r \rangle$, so that its Levi factor is isomorphic to
$\GO(W)\times\GL(r)$. Further let $S_Q=R\cap( \GSp(n)\times Q)$ be
the parabolic subgroup of $R$ whose Levi factor $M_Q$ is isomorphic
to $R_{n, W}\times\GL(r)$, where $R_{n,W}$ is defined in the same
way as $R$, but with respect to $\GSp(n)$ and $\GO(W)$. We denote by
$N_Q$ its unipotent radical. Also let $S_{Q_r}$ be the parabolic
subgroup of $M_Q$ whose Levi factor $M_{Q_r}$ is isomorphic to
$R_{n,W}\times\Gm^r$, \ie $S_{Q_r}=M_Q\cap(\GSp(n)\times Q_r)$. We
denote by $N_{Q_r}$ its unipotent radical. Now let $P_i$ be the
standard parabolic subgroup of $\GSp(n)$ whose Levi factor is
isomorphic to $\Gm^i\times\GSp(n-i)$. Then we define $S_{P_i,Q_r}$
to be the parabolic subgroup of $M_{Q_r}$ whose Levi factor $M_{P_i,
Q_r}$ is isomorphic to $\Gm^i\times R_{n-i,W}\times \Gm^r$, \ie
$S_{P_i,Q_r}=M_{Q_r}\cap(P_i\times Q_r)$. We write a typical element
in $M_{P_i,Q_r}$ by $(\alpha_1,\cdots,\alpha_i, (g,h),
\beta_1,\dots,\beta_r)$. Notice we have the inclusions
$S_{P_i,Q_r}\subset M_{Q_r} \subset S_{Q_r} \subset M_{Q}$. Also we
can set $P_0=\GSp(n)$ and so $S_{P_0,Q_r}=M_{P_0,Q_r}=M_{Q_r}$. Now
the unnormalized Jacquet module of $\omega_{n,X}$ is computed as
follows, which is nothing but Lemma 4 of \cite{HST} with the
notations adjusted to ours.

\begin{prop}\label{P:Jacquet}
The unnormalized Jacquet module $J=J(\omega_{n,X})_{N_Q}$ of
$\omega_{n,X}$ with respect to $N_Q$ has a filtration
\[
    0=J^{(s+1)}\subset J^{(s)}\subset\dots\subset J^{(1)}\subset J^{(0)}=J
\]
of $M_Q$-modules, where $s=\min\{n,r\}$. Let
$I^{(i)}=J^{(i)}/J^{(i+1)}$. Then the unnormalized Jacquet module
$I^{(i)}_{N_{Q_r}}$ of $I^{(i)}$ with respect to $N_{Q_r}$, which is
an $M_{Q_r}$-module, is given by
\[
    {I^{(i)}_{N_{Q_r}}}=\Ind_{S_{P_i,Q_r}}^{M_{Q_r}}\sigma_{i,r},
\]
where $\sigma_{i,r}$ is given by the representation of $M_{P_i,Q_r}$
which is of the form
\[
    \begin{aligned}
    (\alpha_1,\cdots,\alpha_i, (g,h), \beta_1,\dots,\beta_r)\mapsto
    &|\nu(g)|^{nr/2-ni-li/4}|\alpha|^{n+l/2}(\alpha,(-1)^{l/2}D_W)|\beta|^n\\
    &\cdot\prod_{j=1}^{i}\mu_{r-i+j}(\alpha_{i-j+1}^{-1}\beta_{r-i+j}\nu(g))\omega_{n-i,W}(g,h)
    \end{aligned}
\]
for some characters $\mu_{r-i+j}$, where
$\alpha=\alpha_1\cdots\alpha_i$, $\beta=\beta_1\cdots\beta_{r-i}$,
$D_W=\disc W$ and $(,)$ is the Hilbert symbol.
\end{prop}

\begin{rmk}
In the above notations, if one of $n-i$ and $W$ is zero,
$\omega_{n-i,W}$ is taken to be the trivial representation. If $n-i$
is zero, we write a typical element in $M_{P_i,Q_r}$ by
$(\alpha_1,\cdots,\alpha_i, (h), \beta_1,\dots,\beta_r)$ where
$h\in\GO(W)$, and we have to replace $\nu(g)$ by $\nu(h)^{-1}$ in
the above formula. If $W$ is zero, we write a typical element in
$M_{P_i,Q_r}$ by $(\alpha_1,\cdots,\alpha_i, (g),
\beta_1,\dots,\beta_r)$, where $g\in\GSp(n-i)$, and $g$ acts as in
the above formula. If both $n-i$ and $W$ are zero, we have
$M_{P_i,Q_r}\cong\Gm^i\times\Gm\times\Gm^r$ and write a typical
element by $(\alpha_1,\cdots,\alpha_i, (\lambda),
\beta_1,\dots,\beta_r)$ if, for the natural projection
$\iota:M_{P_i,Q_r}\ra P_i$, we have $\nu(\iota(\alpha_1,\cdots,\alpha_i,
(\lambda), \beta_1,\dots,\beta_r))=\lambda$, and we have to replace
$\nu(g)$ by $\lambda$.
\end{rmk}

\begin{rmk}
Although this is a small point, the reader should notice that the
choice of the parabolic $R_{P_i, Q}$ in \cite{HST} is not quite
correct and should be replaced by our $S_{P_i, Q_r}$. Also in
\cite{HST} there is a misprint for the index of $\alpha$ inside
$\nu_{r-i+j}$.
\end{rmk}

The following lemma will be necessary later.
\begin{lemma}\label{L:Rhomo}
Keeping the above notations, let $\mu$ and $\delta$ be admissible
representations of $Q_r$ and $P_i$, respectively. Then the natural
map $\Ind_{P_i\times
Q_r}^{\GSp(n)\times\GO(V)}(\delta\otimes\mu)\rightarrow\Ind_{R\cap(P_i\times
Q_r)}^R(\delta\otimes\mu)$ is an injective $R$-homomorphism.
\end{lemma}
\begin{proof}
Let $F\in\Ind_{P_i\times
Q_r}^{\GSp(n)\times\GO(V)}(\delta\otimes\mu)$ and $\bar{F}$ its
image under the natural map, namely $\bar{F}=F|_R$. Assume
$\bar{F}=0$. Then for each $(g,h)\in\GSp(n)\times\GO(V)$, let us
define
\[
    u=u(g,h)=\begin{pmatrix} I_n & O\\
                           O  & \nu(g)\nu(h)I_n\\
            \end{pmatrix}\in\GSp(n),
\]
where $I_n$ is the $n\times n$ identity matrix. Then $(u,1)\in
P_i\times Q_r$, $\nu(u)=\nu(g)\nu(h)$, and $(u^{-1}g, h)\in R$. So
we have
\[
    F(g,h)=F(u u^{-1}g,h)
    =(\delta\otimes\mu)(u,1)F(u^{-1}g,h)
    =(\delta\otimes\mu)(u,1)\bar{F}(u^{-1}g,h)=0.
\]
Thus the map is injective.
\end{proof}


\subsection{Computation of local parameters}


We will compute the local parameters of unramified theta lifts from
$\GO(X)$ to $\GSp(1)(=\GL(2))$. First assume $d\neq 1$. First notice
the following.
\begin{rmk}\label{R:Galois}
Assume the extension $E/F$ is unramified. Then it is easy to see
that any unramified character $\eta$ on $E^{\times}$ is Galois
invariant and written as $\eta=\chi\circ N$ for some unramified
character $\chi$ on $F^{\times}$. Accordingly if
$\tau=\nInd(\eta,\eta')$ is unramified, then $\eta=\chi_1\circ N$
and $\eta'=\chi_2\circ N$ for unramified characters $\chi_1$ and
$\chi_2$ on $F^{\times}$.
\end{rmk}

Then we have

\begin{prop}\label{P:dnot1}
Assume the extension $E/F$ is unramified, and
$\sigma\in\Irr(\GO(X))$ is unramified, \ie $\sigma =\pi^+$ for some
unramified $\pi =\pi(\tau, \chi)$. (So by Remark \ref{R:Galois}
$\pi=\pi^c$.) If $\sigma$ corresponds to $\Pi\in\Irr(\GSp(1))$ and
$\Pi$ is unramified, then $\tau$ is the base change lift of $\Pi$.
In particular such $\Pi$ is unique. Moreover the central character
of $\Pi$ is $\chi_{E/F}\chi$.
\end{prop}
\begin{proof}
Assume $\sigma$ corresponds to unramified $\Pi$. Since $\tau$ is
unramified, we can write $\tau=\nInd_P^{\GL(2,E)}(\eta\otimes\eta')$
for unramified characters $\eta$ and $\eta'$ on $E^{\times}$, or by
using unnormalized induction,
$\tau=\Ind_P^{\GL(2,E)}(\tilde{\eta}\otimes\tilde{\eta}')$, where
$\tilde{\eta}=|\cdot|_E^{1/2}\eta$ and
$\tilde{\eta}'=|\cdot|_E^{-1/2}\eta'$. (Here note that
$|\cdot|_E=|\cdot|\circ N$.) By Remark \ref{R:Galois},
$\eta=\chi_1\circ N$ and $\eta'=\chi_2\circ N$ for unramified
characters $\chi_1$ and $\chi_2$ on $F^{\times}$ with
$\chi_1\chi_2=\chi$. Then
$\pi=\pi(\chi,\tau)=\Ind_Q^{\GSO(X)}(\mu)$, where $Q$ is the
parabolic preserving the flag
\[
    \langle \begin{pmatrix} 0&\sqrt{d}\\0&0\end{pmatrix} \rangle,
\]
and $\mu$ is defined by
\[
    \mu(\begin{pmatrix} \beta_1&\ast&\ast&\ast\\
                0&a&b\ d&\ast\\
                0&b&a&\ast\\
                0&0&0&\lambda {\beta_1}^{-1}\end{pmatrix})
    =(\tilde{\eta}'/\chi)(\lambda {\beta_1}^{-1})\tilde{\eta}(a+b \sqrt{d})\\
    =(|\cdot|\chi_1/\chi_2)(\beta_1)(|\cdot|^{-1/2}\chi_2)(\lambda),
\]
where the matrix representation of each element of $\GO(X)$ is with
respect to the ordered basis
\[
    \begin{pmatrix} 0&\sqrt{d}\\0&0 \end{pmatrix},
    \begin{pmatrix} 1&0\\0&1 \end{pmatrix},
    \begin{pmatrix} \sqrt{d}&0\\0&-\sqrt{d} \end{pmatrix},
    \begin{pmatrix} 0&0\\(2/d)\sqrt{d}&0 \end{pmatrix}.
\]
(Here we view $\tilde{\eta}'$ as a character on $F^{\times}$ via
$F^{\times}\hookrightarrow E^{\times}$. See also p.279-278 in
\cite{Rob01}.) The second equality can be seen as follows. First
notice that the middle block
\[
    h=\begin{pmatrix} a&b\ d\\
              b&a
      \end{pmatrix}
\]
is identified with $E^{\times}$ by $h\leftrightarrow a+b\sqrt{d}$.
In particular the Levi factor of $Q$ is isomorphic to
$\GO(W)\times\Gm$ where $W$ is the quadratic space $E$ equipped with
the quadratic form $-N$. Thus $\lambda=N(h)$ and so
$\tilde{\eta}(a+b\sqrt{d})=(|\cdot|_E^{1/2}\eta)(h)=(|\cdot|^{1/2}\chi_1)(N(h))
=(|\cdot|^{1/2}\chi_1)(\lambda)$.

Now assume $\Pi$ corresponds to $\sigma$, so there is a non-zero
$R$-homomorphism $\omega_{1,X} \longrightarrow \Pi\otimes\sigma$.
Notice that there is an injective $\GO(X)$-homomorphism
$\sigma\hookrightarrow\Ind_{\GSO(X)}^{\GO(X)}\sigma
(\cong\Ind_Q^{\GO(X)}\mu)$. Then by composing with this injection,
we have a non-zero $R$-homomorphism
\[
    \omega_{1,X} \longrightarrow \Pi\otimes\Ind_Q^{\GSO(X)}\mu.
\]
By Lemma \ref{L:Rhomo}, we have a non-zero $R$-homomorphism
\[
    \omega_{1,X} \longrightarrow \Ind_{R\cap(\GSp(1)\times Q)}^R \Pi\otimes\mu.
\]
Then if we take the Jacquet module of $\omega_{1,X}$ with respect to
$N_Q$, the Frobenius reciprocity together with Proposition
\ref{P:Jacquet} gives $M_Q$-homomorphisms
\[
    0\subset J^{(1)}\subset J^{(0)} \xrightarrow{\quad\varphi\quad}\Pi\otimes\mu.
\]
(Here $\Pi$ is actually $\Pi|_{\GSp(1)^{+}}$, where
$\GSp(1)^{+}=\{g\in\GSp(1): \nu(g)\in\nu(\GO(X))\}$.) It is easy to
see that $\ker\varphi\nsupseteq J^{(1)}$. (Indeed, if
$\ker\varphi\supseteq J^{(1)}$, by applying Proposition
\ref{P:Jacquet}, we will get a contradiction.)

So by restricting $\varphi$ to $J^{(1)}$, Proposition
\ref{P:Jacquet} together with $M_{Q_1}=M_Q$ gives a non-zero
$M_Q$-homomorphism
\[
   J^{(1)}(=I^{(1)})=\Ind_{S_{P_1,Q}}^{M_Q}\sigma_{1,1}\xrightarrow{\quad\varphi'\quad}\Pi\otimes\mu.
\]
By considering the action of the elements of the form
$(1,(1),\beta_1)$ on $\Ind_{S_{P_1,Q}}^{M_Q}\sigma_{1,1}$ and
$\Pi\otimes\mu$, we see that
\[
    \mu_1=\chi/\tilde{\eta}'=|\cdot|\chi_1/\chi_2,
\]
where $\mu_1$ is as in Proposition \ref{P:Jacquet}.

Assume $\Pi=\Ind_{P_1}^{\GSp(1)}\delta$ for the character $\delta$
on $P_1$ given by
\[
    \delta(\begin{pmatrix} \alpha_1&\ast\\
                    0&\lambda\alpha_1^{-1}
    \end{pmatrix})
    =\delta_1(\alpha_1)\delta_2(\lambda).
\]
Then $\Pi\otimes\mu\cong\Ind_{P_1\times Q}^{\GSp(1)\times
Q}(\delta\otimes\mu)$. By Lemma \ref{L:Rhomo} $\varphi'$ gives rise
to a non-zero $M_Q$-homomorphism
\[
    \Ind_{S_{P_1,Q}}^{M_Q}\sigma_{1,1}\xrightarrow
    {\quad\varphi'\quad}\Ind_{S_{P_1,Q}}^{M_Q}(\delta\otimes\mu),
\]
which we also call $\varphi'$. Let us view $\varphi'$ as an
$R$-homomorphism by restricting to the elements of the form
$((g,h),1)$, and so we have a non-zero $R$-homomorphism
\[
   \Ind_{R_{1,W}}^R\sigma_{1,1}\xrightarrow
   {\quad\varphi'\quad}\Ind_{R_{1,W}}^R(\delta\otimes\mu),
\]
where $\sigma_{1,1}$ and $\delta\otimes\mu$ are restricted to
$R_{1,W}$. This $\varphi'$ can be made into a non-zero
$\GSp(1)$-homomorphism $\varphi''$ via the diagram
\[
    \xymatrix{ \Ind_{R_{1,W}}^R\sigma_{1,1}\ \ar[d]^{\varphi'}\ar@{^{(}->}[r]^{\tilde{\Phi}}
         & \Ind_{P_1^{+}}^{\GSp(1)^{+}} \widehat{\sigma_{1,1}}\ \ar[d] \ar@{^{(}->}[r]
         & \Ind_{P_1}^{\GSp(1)} \widehat{\sigma_{1,1}} \ar[d]^{\varphi''}\\
           \Ind_{R_{1,W}}^R(\delta\otimes\mu)\ \ar@{^{(}->}[r]
         & \Ind_{P_1^{+}}^{\GSp(1)^{+}}(\widehat{\delta\otimes\mu})\ \ar@{^{(}->}[r]
         & \Ind_{P_1}^{\GSp(1)}(\widehat{\delta\otimes\mu}),
         }
\]
where all the horizontal arrows are injective and the vertical
arrows are non-zero intertwining operators for the corresponding
groups. This diagram is given as follows. First, define
$\iota:R\ra\GSp(1)^{+}$ by $(g,h)\mapsto g$. Also let
$\widehat{\sigma_{1,1}}$ be the character on $P_1^{+}(=\GSp(1)\cap
P_1)$ defined by
\[
    \begin{pmatrix}\alpha_1&\ast\\0&\lambda\alpha_1^{-1}\end{pmatrix}
    \mapsto\sigma_{1,1}(\alpha_1,(h),1)
    =|\nu(h)||\alpha_1|^2(\alpha_1,-D_W)\mu(\alpha_1^{-1}\nu(h)^{-1}),
\]
where $h\in\GO(W)$ is such that $\nu(h)=\lambda^{-1}$. This is
well-defined. (Also note that $(\alpha_1,-D_W)=(\alpha_1,d)$ is the
quadratic character for the quadratic extension $E/F$ and we put
$(\alpha_1,d)=\chi_{E/F}(\alpha_1)$). This gives us a natural map
$\Phi:\sigma_{1,1}\ra\widehat{\sigma_{1,1}}$ that respects the
actions of $R_{1,W}$ and $P_1^{+}$ via $\iota$. Then $\Phi$ can be
extended to a non-zero map
\[
   \tilde{\Phi}:\Ind_{R_{1,W}}^R\sigma_{1,1}\longrightarrow\Ind_{P_1}^R\widehat{\sigma_{1,1}}
\]
defined by $\tilde{\Phi}(F)(g)=\Phi(F((g,h),1))$, where $h\in\GO(W)$
is such that $\iota(g,h)=g$. This is well-defined, injective, and
respects the actions of $R$ and $\GSp(1)^{+}$ via $\iota$. We can
similarly define a non-zero injective map
\[
   \Ind_{R_{1,W}}^R(\delta\otimes\mu)\longrightarrow
   \Ind_{P_1^{+}}^{\GSp(1)^{+}}(\widehat{\delta\otimes\mu}).
\]
Here $\widehat{\delta\otimes\mu}(\left(\begin{smallmatrix}\alpha_1 &
\ast\\0 &
\lambda{\alpha_1}^{-1}\end{smallmatrix}\right))=\delta_1(\alpha_1)\delta_2(\lambda)
(|\cdot|^{-1/2}\chi_2)(\lambda^{-1})$. To see the left square of the
diagram, notice that, in general, if $\gamma$ is a character on
$P_1$, there is an injective map
$\Ind_{{P_1}^{+}}^{\GSp(1)^{+}}\gamma\ra\Ind_{P_1}^{\GSp(1)}\gamma$
given by $F\mapsto\tilde{F}$ where
$\tilde{F}(g)=\gamma(\left(\begin{smallmatrix}1&0\\0&\nu(g)^{-1}\end{smallmatrix}\right))
F((\left(\begin{smallmatrix}1&0\\0&\nu(g)\end{smallmatrix}\right)g)$.
This map respects the actions of $\GSp(1)^{+}$ and $\GSp(1)$ via the
inclusion $\GSp(1)^{+}\hookrightarrow\GSp(1)$. Since both
$\widehat{\sigma_{1,1}}$ and $\widehat{\delta\otimes\mu}$ can be
viewed as characters on $P_1$ in the obvious way, we have the
injective maps as in the above diagram.

Now let us switch to the notation $P=\left(\begin{smallmatrix}
a&\ast\\0&d\end{smallmatrix}\right)$ for the parabolic. Then
\[
    \widehat{\sigma_{1,1}}(\begin{pmatrix}a&\ast\\0&d\end{pmatrix})
    =(|\cdot|\chi_{E/F})(a)(\chi_1/\chi_2)(d)
\]
and
\[
    \widehat{\delta\otimes\mu}(\begin{pmatrix}a&\ast\\0&d\end{pmatrix})
    =(\delta_1\delta_2|\cdot|^{1/2}/\chi_2)(a)(|\cdot|^{1/2}\delta_2/\chi_2)(d).
\]
Note that, by twisting by $\chi_2/|\cdot|^{1/2}$, $\varphi''$ gives
rise to a non-zero $\GL(2)$-homomorphism
\[
   (\chi_2/|\cdot|^{1/2})\otimes\Ind_P^{\GL(2)}(\widehat{\sigma_{1,1}}) \longrightarrow
   (\chi_2/|\cdot|^{1/2})\otimes\Ind_P^{\GL(2)}(\widehat{\delta\otimes\mu}),
\]
which can be rewritten as
\[
    \nInd_P^{\GL(2)}(\chi_{E/F}\chi_2\otimes\chi_1)\longrightarrow\Pi
\]
Therefore $\Pi\cong\nInd_P^{\GL(2)}(\chi_{E/F}\chi_2\otimes\chi_1)$,
and so $\tau$ is the base change lift of $\Pi$ whose central
character is $\chi_{E/F}\chi_1\chi_2=\chi_{E/F}\chi$.
\end{proof}

The case $d=1$ is well known. Namely,

\begin{prop}\label{P:d1}
Assume $d=1$ and $\sigma\in\Irr(\GO(X))$ is unramified, \ie $\sigma
=\pi^+$ for some unramified $\pi =\pi(\tau_1, \tau_2)$ with
$\tau_i$'s unramified representations of $\GL(2,F)$. If $\sigma$
corresponds to $\Pi\in\Irr(\GSp(1))$, then
$\tau_1\cong\tau_2\cong\Pi$. In particular, such $\Pi$ is unique.
\end{prop}

This proposition is known to be true even for the the ramified case.
But it might be a good exercise to prove the unramified case in the
way we did for the $d\neq 1$ case.


\section{Proof of Theorem  \ref{T:main1} and \ref{T:main2}}\label{S:Main}


Now we are ready to prove our main theorems. First we give a proof of
Theorem \ref{T:main1}, the major arguments of which are essentially
given in \cite{Rob01}.

\begin{proof}[Proof of Theorem \ref{T:main1}]
Let $\sigma$ be a cuspidal automorphic representation of
$\GO(X,\A_F)$ satisfying the hypotheses of Theorem \ref{T:main1}. Also let us assume that $\sigma$ is extended from an irreducible cuspidal automorphic representation
$\pi$ of $\GSO(X,\A_F)$, \ie $\sigma=(\pi,\delta)$ for some $\delta$
as discussed at the end of Section \ref{S:GSO}. Also let $\sigma_1$ be an irreducible
constituent of the space $\{f|_{\OO(X,\A_F)} : f\in V_{\sigma}\}$
viewed as an automorphic representation of $\OO(X,\A_F)$.

First let us show the if part. So assume
$\sigma\cong\otimes\sigma_v$ has the property that $\sigma_v$ has a
local theta lift to $\GSp(2,F_v)$ for all $v$. Then as discussed in
\cite[Lemma 4.2]{Rob96}, $\sigma_1\cong\otimes{\sigma_1}_v$ has the
property that ${\sigma_1}_v$ has a local theta lift to $\Sp(2,
F_v)$. Also note that
\[
    L^{S}(s,\sigma_1)
    =\begin{cases}L^{S}(s, \sigma_1)=L^{S}(s, \tau_1^{\JL}\times{\tau_2^{\JL}}^{\vee})
    \text{ if $d=1$ and $\pi=\pi(\tau_1,\tau_2)$}\\
    L^{S}(s, \tau^{\JL}, \chi^{-1}, \text{Asai})
    \text{ if $d\neq 1$ and $\pi=\pi(\tau,\chi)$}\end{cases},
\]
where $^{\JL}$ indicated the Jacquet-Langlands lift, $L^{S}(s,
\tau_1^{\JL}\times\tau_2^{\JL})$ is the (incomplete) Rankin-Selberg
$L$-function, and $L^{S}(s, \tau^{\JL}, \chi^{-1}, \text{Asai})$ is
the (incomplete) Asai $L$-function. (See Lemma 8.1 of \cite{Rob01}
for the proof.) It is well-known that the Rankin-Selberg
$L$-function does not vanish at $s=1, 2$. The Asai $L$-function is
also non-vanishing at $s=1, 2$. This can be shown as follows. First,
the non-vanishing at $s=1$ is given in, say, \cite{Flicker88}, which
is also discussed in \cite[p.302]{Rob01}. The non-vanishing at $s=2$
follows from absolute convergence of the Euler product for
$\Re(s)>1+\frac{2}{9}$. The proof for this is elementary. Namely, by
the explicit description of the unramified factor of the Euler
product of the Asai $L$-function found in, say,
\cite[p.200]{Flicker93}, and the known estimate of the bound of the
Langlands parameters given by \cite{Kim02}, one sees that each
unramified Euler factor of the Asai $L$-function is the product of
four factors of the form
\[
    (1-a_v q_v^{-s})^{-1} \quad \text{with} \quad |a_v|<q_v^{2/9}.
\]
Then the convergence of this product can be easily proven by using
Lemma 2 in p.187 and the first paragraph of p.188 of \cite{Marcus}.
Thus by Theorem \ref{T:main1} and Lemma \ref{L:comparison}, we see
that the global theta lift of $\sigma$ to $\GSp(2, \A_F)$ does not
vanish.

Conversely, assume the global theta lift $\Theta_2(V_\sigma)$ to $\GSp(2,\A_F)$ does not vanish. If $\Theta_2(V_\sigma)$ is not contained in the space of cusp forms, then by the persistence principle the global theta lift $\Theta_1(V_\sigma)$ to $\GSp(1,\A_F)$ is non-zero, and is contained in the space of cusp forms. Now let $(\Pi, V_\Pi)$ be an irreducible constituent of the lift $\Theta_1(V_\sigma)$. Then each $\sigma_v$ and $\Pi_v^\vee$ correspond under the local theta correspondence for the similitude in the sense defined at the beginning of Section \ref{S:similitude}. Hence by the persistence principle, for each $v$, $\theta_2(\sigma_v)\neq 0$. If $\Theta_2(V_\sigma)$ is contained in the space of cusp forms, similarly one sees that $\theta_2(\sigma_v)\neq 0$ for each $v$.
\end{proof}

Next we prove Theorem \ref{T:main2}. For this, we need

\begin{prop}\label{P:Asai}
Let $E/F$ be a quadratic extension of global fields, and $\tau$ an irreducible
cuspidal automorphic representation of $\GL(2, \A_E)$ whose central
character is $\chi\circ N^E_F$ for a Hecke character $\chi$ on
$\A_F^{\times}$. Then the incomplete Asai $L$-function $L^S(s, \tau,
\chi^{-1}, \text{Asai})$ has a pole at $s=1$ if and only if $\tau$
is the base change lift of a cuspidal automorphic representation
$\tau_0$ of $\GL(2,\A_F)$ whose central character is
$\chi_{E/F}\chi$.
\end{prop}
\begin{proof}
The case where $\chi$ is trivial is essentially done in P.311 of
\cite{Flicker88}. (Also see \cite{Flicker95} for the archimedean
assumption imposed on \cite{Flicker88}.) We can extend his method to
the case where $\chi$ is non-trivial by making the following small
modifications to \cite{Flicker88}, though we restrict to the case
$G=\GL(2)$. (Indeed, if $G\neq\GL(2)$, our argument would not work.)
So in the following, all the notations are as in \cite{Flicker88}.

First for each $\phi\in V_{\tau}$, let us define
\[
    I(s,\Phi,\phi,\chi)=\int_{Z(\A_F)G(F)\backslash G(\A_F)}
    E(g,\Phi,s)\phi(g)\chi^{-1}(\det g) dg,
\]
where $\Phi$ is a Schwartz-Bruhat function on $\A_F^2$ and
$E(g,\Phi,s)$ is the Eisenstein series as in \cite{Flicker88}.
(Compare this with the integral in \cite[p.302]{Flicker88}.) Then as
in \cite{Flicker88}, we see that $I(s,\Phi,\phi,\chi)$ has a pole at
$s=1$ if and only if $\hat{\Phi}(0)\neq0$ and $\int
\phi(g)\chi^{-1}(\det g)dg\neq 0$.

Next consider the integral
\[
    \Psi(s,\Phi,W,\chi)=\int_{N(\A_F)\backslash G(\A_F)}
    W(g)\Phi(\epsilon g)|g|^s\chi^{-1}(\det g) dg,
\]
where $W$ is the Whittaker vector as defined in
\cite[p.302]{Flicker88}. Then by the same computation as in
\cite[p.303]{Flicker88}, we have
\[
    I(s,\Phi,\phi, \chi)=\Psi(s,\Phi,W,\chi),
\]
and
\[
    \Psi(s, \Phi, W, \chi)=A(s, \Phi_S, W_S, \chi_S)L^S(s, \tau, \chi^{-1},
    \text{Asai}).
\]
for some function $A(s, \Phi_S, W_S, \chi_S)$ with
$\hat{\Phi}(0)\neq0$. (See \cite[p.310]{Flicker88}.) Therefore the
(incomplete) Asai $L$-function $L^S(s, \tau, \chi^{-1},\text{Asai})$
has a pole at $s=1$ if and only if $\Psi(s, \Phi, W, \chi)$ has a
pole at $s=1$, \ie if and only if $\int \phi(g)\chi^{-1}(\det
g)dg\neq 0$. The rest of the proof is basically the same as
\cite{Flicker88}, and is left to the reader.
\end{proof}

Now we have

\begin{proof}[Proof Theorem \ref{T:main2}]
Assume $\sigma$ is an irreducible cuspidal automorphic representation of
$\GO(X, \A_F)$ with the global theta lift $\Theta_2(V_{\sigma})$
non-vanishing. (We write $\sigma=(\pi,\delta)$ as in Section
\ref{S:GSO}.) Also let $\sigma_1$ be as in the proof of Theorem
\ref{T:main1}. Then by Lemma \ref{L:comparison} the global theta
lift $\Theta_2(V_{\sigma_1})$ to $\Sp(2,\A_F)$ does not vanish. Now
assume $\Theta_1(V_{\sigma_1})=0$ so that $\Theta_2(V_{\sigma_1})$
provides an irreducible cuspidal automorphic representation of $\Sp(2,\A_F)$,
which we denote by $\Pi$. Then by functoriality of the unramified
theta correspondence, we have
$L^S(\Pi,s,\chi_X)=\zeta^S_F(s)L^S(\sigma_1, s)$. It is well known
that $\zeta^S_F(s)$ has a simple pole at $s=1$. So if $L^S(\sigma_1,
s)$ has a simple pole at $s=1$, then $L^S(\Pi,s,\chi_X)$ would have
a double pole, which contradicts to the fact that the poles of the
standard (incomplete) Langlands $L$-function of $\Sp(n)$ are at most
simple. (See Remark after Proposition 1.7 of \cite{Ikeda92}, and
also Theorem 7.2.5 of \cite{Kudla-Rallis94}.) Therefore if
$\Theta_1(V_{\sigma_1})=0$, then $L^S(\sigma_1, s)$ can not have a
pole at $s=1$. Hence if $d=1$ and so $\pi=\pi(\tau_1,\tau_2)$, then
$L^{S}(s, \tau_1^{\JL}\times{\tau_2^{\JL}}^{\vee})$ has a pole at
$s=1$ \ie $\tau_1=\tau_2$. Also if $d\neq 1$ and so
$\pi=\pi(\tau,\chi)$, then $L^S(s, \tau^{\JL},\chi^{-1},
\text{Asai})$ has a pole at $s=1$, and thus by Proposition
\ref{P:Asai} we see that $\tau^{\JL}$ is the base change lift of a
cuspidal automorphic representation $\tau_0$ of $\GL(2,\A_F)$ whose
central character is $\chi_{E/F}\chi$.

Conversely, assume $\Theta_1(V_{\sigma})\neq0$, and so it gives rise
to a cuspidal automorphic representation $\Pi$ of $\GL(2,\A_F)$.
First assume $d=1$ and so $\pi=\pi(\tau_1,\tau_2)$. Then note that
for almost all $v$, we have $\pi_v=\pi({\tau_1}_v, {\tau_2}_v)$
where ${\tau_1}_v$ and ${\tau_2}_v$ are spherical representations of
$\GL(2,F_v)$. Then by Proposition \ref{P:d1} we have
${\tau_1}_v\cong{\tau_2}_v$. So by strong multiplicity one theorem,
$\tau_1\cong\tau_2$. Also in this case we know that
$(\Pi_v)^{\vee}\cong{\tau_1}_v\cong{\tau_2}_v$ for almost all $v$.
Thus we have $\Pi^{\vee}\cong\tau_1\cong\tau_2$. Moreover since
$\Theta_1(V_{\sigma})$ is closed under the right action of
$\GL(2,\A_F)$, multiplicity one theorem gives us
$\Theta_1(V_{\sigma})=V_{\Pi}=V_{\pi_1^{\vee}}=V_{\pi_2^{\vee}}$. We
can similarly prove the $d\neq 1$ case by using strong multiplicity
one and Proposition \ref{P:d1} and \ref{P:dnot1}. This completes the
proof.
\end{proof}


\subsection{A remark for the $d\neq1$ case}


Our theorem for the $d\neq1$ case reveals the following interesting
phenomenon. Let $\tau$ be an irreducible cuspidal automorphic representation of
$\GL(2,\A_E)$ that is the base change lift of an irreducible cuspidal automorphic
representation $\tau_0$ of $\GL(2,\A_F)$. (Note that such $\tau_0$
is unique.) Assume the central character of $\tau_0$ is $\chi$, and
so the central character of $\tau$ is $\chi\circ N_F^E$. Then we can
identify $\tau$ with an irreducible cuspidal automorphic representation
$\pi=\pi(\tau,\chi)$ of $\GSO(X,\A_F)$ whose central character is
$\chi$. Since $\tau$ is generic, $\pi$ can be always extended to an irreducible
cuspidal automorphic representation $\sigma=(\pi,\delta)$ of
$\GO(X,\A_F)$ such that $\Theta_2(V_{\sigma})\neq0$. (See Appendix
\ref{S:disconnected}.) In this case by Theorem \ref{T:main2} we have
$\Theta_1(V_{\sigma})=0$. Thus $\Theta_2(V_{\sigma})$ is cuspidal.

However since $\chi\circ N=\chi_{E/F}\chi\circ N_F^E$, we can also
identify $\tau$ with an irreducible cuspidal automorphic representation
$\pi'=\pi'(\tau,\chi_{E/F}\chi)$ of $\GSO(X,\A_F)$ whose central
character is $\chi_{E/F}\chi$. Again we can extend $\pi'$ to
$\sigma'=(\pi', \delta')$ such that $\Theta_2(V_{\sigma'})\neq0$.
Then this time, by Theorem \ref{T:main2} we have
$\Theta_1(V_{\sigma'})=V_{\tau_0^{\vee}}\neq0$, and
$\Theta_2(V_{\sigma'})$ is non-zero but non-cuspidal. This
observation can be summarized as follows.
\[
    \xymatrix@R=2pt{
            &&\pi=\pi(\chi,\tau)\ar[r]&\sigma=(\pi,\delta)\ar[r]^-{\txt\scriptsize{theta\\lift}}&
                                *{\begin{cases} \Theta_2(V_{\sigma})\neq0\\
                                                \Theta_1(V_{\sigma})=0
                                                \end{cases}}\\
        \tau_0\ar[r]^-{\txt\scriptsize{base\\change}}&\tau\ar[ur]\ar[dr]&&&\\
            &&\pi'=\pi'(\chi_{E/F}\chi,\tau)\ar[r]&\sigma'=(\pi',\delta')\ar[r]^-{\txt\scriptsize{theta\\lift}}&
                                *{\begin{cases} \Theta_2(V_{\sigma'})\neq0\\
                                                \Theta_1(V_{\sigma'})=V_{\tau_0^{\vee}} \end{cases}}
    }
\]
(Of course this diagram makes sense only with the assumption that
$\tau$ is cuspidal. It is known that $\tau$ is cuspidal if and only
if $\tau_0$ is the automorphic induction of a Hecke character of
$E$.)


\appendix
\section{Proof of Theorem 1.4}\label{S:disconnected}


In this appendix, we will prove Theorem \ref{T:main3}. Let us start
with the following theorem due to \cite{H-PS83}.

\begin{prop}\label{P:Howe_PS}
Assume $\sigma$ is a generic irreducible cuspidal automorphic representation of
$\OO(X, \A_F)$ or $\GO(X, \A_F)$. Also let
$T=\OO(X)\slash\SO(X)\cong\GO(X)\slash\GSO(X)$. Then there is a
unitary character $\mu$ of $T(F)\backslash T(\A_F)$ such that the
space $\mu\otimes V_{\sigma}=\{\mu f : f\in V_{\sigma}\}$ has a
non-zero theta lift $\Theta_2(\mu\otimes V_{\sigma})$ to
$\Sp(2,\A_F)$ or $\GSp(2,\A_F)$, respectively.
\end{prop}
\begin{proof}
The isometry case is essentially Theorem 8.1 of \cite{H-PS83}, which
follows from Theorem 5.7 of the same paper. However, we should
mention that in \cite{H-PS83} they state that $\mu$ is on $T(\A_F)$,
but from the proof of Theorem 5.7 it is easy to see that $\mu$ is
actually on $T(F)\backslash T(\A)$. Indeed if $\mu$ is not trivial
on $T(F)$, $\mu f$ can not be automorphic. The similitude case
follows immediate from the isometry case and Lemma
\ref{L:comparison}.
\end{proof}

We should also mention the following, which is due to Roberts
\cite{Rob99, Rob01}.

\begin{prop}\label{P:Rob_local}
Let $\pi_v$ be an irreducible admissible representation of $\GSO(X,F_v)$. Then
at least one of $\pi_v^{+}$ and $\pi_v^{-}$ has non-zero theta lift
to $\GSp(2,F_v)$. Moreover, if $\pi_v$ is spherical, then
$\pi_v^{+}$ always has non-zero theta lift.
\end{prop}
\begin{proof}
The first part is one of the main theorems of \cite{Rob99}. (Note
that the proof works even if $v|\infty$.) The second part is proven
in \cite[Proposition 4.3]{Rob01}.
\end{proof}

Finally we can prove Theorem \ref{T:main3}

\begin{proof}[Proof of Theorem \ref{T:main3}]
(1). First consider any extension $\sigma=(\pi,\delta)$ of $\pi$.
Then by Proposition \ref{P:Howe_PS}, there exists a unitary
character $\mu$ so that $\mu\otimes V_{\sigma}$ has non-zero theta
lift. Now for each $f\in V_{\sigma}$, clearly $f|_{\GSO(X,
\A_F)}=(\mu f)|_{\GSO(X, \A_F)}$, which means that
$V_{\sigma}^{\circ}=(\mu\otimes V_{\sigma})^{\circ}$, \ie both
$V_{\sigma}$ and $\mu\otimes V_{\sigma}$ lie above $V_{\pi}$. Hence
there is an extension index $\delta'$ of $\pi$ so that $\mu\otimes
V_{\sigma}=V_{\sigma'}$ for $\sigma'=(\pi, \delta')$, which we
re-choose to be $\sigma$.

(2). This immediately follows from Proposition
\ref{P:gobal_disconnected}, Theorem \ref{T:main1}, and Proposition
\ref{P:Rob_local}.
\end{proof}

\end{document}